\documentclass[a4paper,oneside,11pt]{article}
\usepackage{graphicx,afterpage}
\usepackage{bm}\usepackage{soul}
\usepackage{empheq}\usepackage{mathtools}
\usepackage[top=0.6in, bottom=1in, left=1in, right=0.8in]{geometry}
\usepackage{mathrsfs,color}
\usepackage{amsfonts}\usepackage{multirow}
\usepackage{amssymb}\usepackage{caption,comment}
\usepackage{amsmath}\usepackage{float}
\usepackage{amssymb,amsthm}
\usepackage{graphicx,afterpage,cancel}

\usepackage{epstopdf}
\usepackage[pdftex,colorlinks=true]{hyperref}
\hypersetup{CJKbookmarks,%
bookmarksnumbered,%
colorlinks,%
linkcolor=black,%
citecolor=black,%
plainpages=false,%
pdfstartview=FitH}
\numberwithin{equation}{section}
\numberwithin{figure}{section}

\newcommand\tabcaption{\def\@captype{table}\caption}

\usepackage[numbers]{natbib}
\setcitestyle{authoryear,round}

\definecolor{orange}{RGB}{255,127,0}

\def\d{{\, \rm d}}

\afterpage{\clearpage}
\usepackage{authblk}
\title{A Physics-Informed Auto-Learning Framework for Developing Stochastic Conceptual Models for ENSO Diversity}
\author[1]{Yinling Zhang}
\author[1,*]{Nan Chen}
\author[2]{J\'er\^ome Vialard}
\author[3]{Xianghui Fang}
\affil[1]{Department of Mathematics, University of Wisconsin-Madison, Madison, USA}
\affil[2]{LOCEAN-IPSL, CNRS/IRD/MNHN/Sorbonne Universit\'e, Paris, France}
\affil[3]{Department of Atmospheric and Oceanic Sciences and Institute of Atmospheric Sciences, Fudan University, Shanghai, China}
\affil[*]{Corresponding author: Nan Chen, chennan@math.wisc.edu}
\date{\today}
\begin{document}
\maketitle %\tableofcontents
\abstract{Understanding ENSO dynamics has tremendously improved over the past decades. However, one aspect still poorly understood or represented in conceptual models is the ENSO diversity in spatial pattern, peak intensity, and temporal evolution. In this paper, a physics-informed auto-learning framework is developed to derive ENSO stochastic conceptual models with varying degrees of freedom. The framework is computationally efficient and easy to apply. Once the state vector of the target model is set, causal inference is exploited to build the right-hand side of the equations based on a mathematical function library. Fundamentally different from standard nonlinear regression, the auto-learning framework provides a parsimonious model by retaining only terms that improve the dynamical consistency with observations. It can also identify crucial latent variables and provide physical explanations. Exploiting a realistic six-dimensional reference recharge oscillator-based ENSO model, a hierarchy of three- to six-dimensional models is derived using the auto-learning framework and is systematically validated by a unified set of validation criteria assessing the dynamical and statistical features of the ENSO diversity. It is shown that the minimum model characterizing ENSO diversity is four-dimensional, with three interannual variables describing the western Pacific thermocline depth, the eastern and central Pacific sea surface temperatures (SSTs), and one intraseasonal variable for westerly wind events. Without the intraseasonal variable, the resulting three-dimensional model underestimates extreme events and is too regular. The limited number of weak nonlinearities in the model are essential in reproducing the observed extreme El Ni\~nos and nonlinear relationship between the eastern and western Pacific SSTs.
}

\section{Introduction}
\bibliographystyle{unsrtnat}
El Ni\~no-Southern Oscillation (ENSO) dominates the planet's interannual variability with significant impacts on global weather and climate through atmospheric teleconnections \citep{philander1983nino, ropelewski1987global, klein1999remote, mcphaden2006enso, dai2000global}. In the traditional viewpoint, ENSO consists of alternating periods of anomalously warm El Ni\~no conditions and cold La Ni\~na conditions that peak in the equatorial eastern Pacific (EP). However, in recent decades, many El Ni\~no events have been observed to peak in the central Pacific (CP) \citep{ashok2007nino, kao2009contrasting, kim2012statistical}. The discovery of this new type of ENSO events leads to the El Ni\~no diversity concept \citep{capotondi2015understanding}, which underlines that El Ni\~no events vary continuously between the CP and EP ideal types \citep{larkin2005global, yu2007decadal, ashok2007nino, kao2009contrasting, kug2009two, jin2022toward}. The shift of the warming center can cause significant differences in the air-sea coupling over the equatorial Pacific, changing how ENSO affects the global climate \citep{chen2008nino, jin2008current, barnston2012skill, hu2012analysis, zheng2014asymmetry, fang2015cloud, sohn2016strength, santoso2019dynamics, stuecker2013combination}. In addition to these two major categories, individual ENSO events further exhibit diverse characteristics in spatial pattern, peak intensity, and temporal evolution. This is known as the ENSO complexity \citep{timmermann2018nino}.

Conceptual models of ENSO are useful due to their low dimensionality, which makes them physically explainable, and mathematically tractable, then further provide the test beds for examining different physical hypotheses and interdependence between state variables \citep{jin1997equatorial, jin1997equatorial2}. Several low-order conceptual models have been developed, including the recharge-discharge oscillator \citep{jin1997equatorial}, the delayed oscillator \citep{suarez1988delayed, battisti1989interannual}, the western-Pacific oscillator \citep{weisberg1997western}, and the advective-reflective oscillator \citep{picaut1997advective}. Later, a unified ENSO oscillator motivated by the dynamics and thermodynamics of Zebiak and Cane's coupled ocean-atmosphere model has been built \citep{wang2001unified}. Some other conceptual models have also been constructed to emphasize the nonlinear dynamics of ENSO and explore the ENSO features in the decadal time scale \citep{timmermann2003nonlinear, roberts2016mixed, timmermann2002nonlinear}. These models were mainly proposed based on physical intuition and highlighted one or two specific dynamical features of the ENSO as building blocks. They primarily focused on EP El Ni\~nos, and led to many successes in understanding and predicting their key characteristics.

Comparatively, few models have been built to describe CP events or capture ENSO diversity and complexity. It was suggested in an early work \citep{ren2013recharge} that EP and CP dynamics can both be treated as recharge-discharge oscillators but with different strengths of the thermocline and the zonal advective feedback. Along this direction, the conceptual model in \citep{geng2020two} generalized the original discharge-recharge oscillation for the EP El Ni\~no \citep{jin1997equatorial} by adding a prognostic equation for the CP SST. Leveraging the recharge-discharge oscillator concept, \citep{thual2023enso} proposed a low-dimensional model representing ENSO diversity by replacing the fixed SST index with a zonally variable warm-pool edge index. On the other hand, \citep{fang2018three} extended the deterministic two-region model to a three-region one, including the western Pacific (WP) heat content, CP SST and zonal current, as well as EP SST. Recently, \citep{chen2022multiscale} extended Fang and Mu model \citep{fang2018three} by including seasonality and two stochastic equations representing the effect of state-dependent synoptic atmosphere forcing such as \citep{jin2007ensemble} and decadal fluctuations in the strength of the Walker circulation. The decadal variability modulates the amplitude of the zonal advection in the CP region. It leads to the system alternating between EP-dominant and CP-dominant periods as observed \citep{capotondi2015understanding}. This stochastic multiscale model reproduces many important dynamical and statistical properties of the observed ENSO complexity, including occurrence frequency \citep{chen2022multiscale}. %Note that in this study, this model will be treated as the reference to represent the nature and utilized to testify the feasibility of this approach.

While recent modeling efforts have advanced the preliminary understanding of the diverse features of ENSO, fundamental issues remain in further conceptualizing ENSO diversity and complexity. First, obtaining the minimal model reproducing key observed features of ENSO is crucial to understanding the most fundamental mechanisms contributing to ENSO diversity. Despite its importance, such a model is still lacking. Second, instead of relying entirely on human knowledge for model development, building an automatic machine-based tool that systematically derives appropriate models satisfying specific requirements is of practical significance. Third, by quantifying model error and missing information in an existing model, it is vital to systematically supplement an effective minimal component that facilitates improving the model performance.

To address these gaps, a physics-informed auto-learning framework is developed in this paper. This framework can automatically derive stochastic conceptual models with different complexities that focus on capturing ENSO diversity. Compared to traditional manual model building, the auto-learning approach enables efficient exploration of model structures without being limited by human biases. Specifically, the model structures are determined using rigorous causal inference that discovers the underlying physics from a comprehensive inference exploiting a balance between the observational data and prior knowledge. The physics-informed auto-learning framework is applied to derive a hierarchy of ENSO stochastic conceptual models with different state variables and degrees of freedom focusing on ENSO diversity. These models will be systematically explored and compared using a unified set of validation criteria metrics that assess how essential ENSO properties, including diversity, are reproduced. Several issues will be studied to understand processes crucial for ENSO diversity. First, the framework will be utilized to test model hierarchies to identify the minimal sufficient model that captures the essential dynamical and statistical features of the ENSO complexity. Second, the roles of linear and nonlinear dynamics in describing the ENSO diversity \citep{chen2017simple, an2004nonlinearity, timmermann2003decadal} will be explored. Such a study aims to discover the fundamental dynamics contributing to the ENSO diversity in different regions. Third, the necessity of representing interactions across timescales \citep{chen2022multiscale, jin2020simple}, including intraseasonal westerly wind (WWB) events, will be analyzed. Lastly, the framework enables identifying effective latent variables referring to additional undetermined variables constructed to improve model performance when the models with predetermined state variables show deficiencies. The time series and resultant properties of these latent variables guide identifying suitable missing physical processes to determine the minimal set of variables needed to capture key ENSO features.

This paper is organized as follows. Section \ref{Sec:Method} describes the physics-informed auto-learning framework for generating stochastic conceptual models of ENSO. Section \ref{Sec:Observations} describes the observational data sets, the reference model whose synthetic data will be used in the auto-learning procedure, and the definitions of various types of ENSO. Section \ref{Sec:Metrics} focuses on describing the validation metrics.  The hierarchy of models derived from the auto-learning framework and their intercomparisons are presented in Section \ref{Sec:Model}. Section \ref{Sec:Discussion} exploits the set of models in understanding the mechanisms and exploring the minimum model of ENSO diversity. The paper is concluded in Section \ref{Sec:Conclusion}.

\section{The Physics-Informed Auto-Learning Framework}\label{Sec:Method}
\subsection{Overview}
The physics-informed auto-learning framework aims to extract the known physical knowledge from existing models with the help of additional data to develop new models automatically. It determines the model structure based on causal inference and provides a physics-informed parsimonious model.
An overview of the framework is as follows.
\begin{enumerate}
  \item [Step 1.] An appropriate existing model is needed to provide a long time series for the auto-learning framework. This model is pre-determined by physical knowledge and is calibrated by observations. The model simulation is expected to meet as many desirable validation criteria as possible. In contrast to short observational data, the long model simulation provides robust learning outcomes.
  \item [Step 2.] The quantities of interests, namely the state variables, are pre-determined to derive a new model using the framework.
  \item [Step 3.] A library of candidate functions is developed as a prerequisite for the auto-learning framework. These candidate functions are potential terms appearing on the right-hand side of one or a few equations of the target model. Physical knowledge can be used to assist the development of such a candidate set. %The function library can also be built by including a large number of functions with arbitrary combinations of state variables that do not refer to any physical intuition.
  \item [Step 4.] A causal inference technique is utilized to reveal the underlying relationship between the time evolution of the state variables and candidate functions. It determines the right-hand side of the model equations and leads to a physically explainable parsimonious model.
  \item [Step 5.] A simple maximum likelihood method is adopted to estimate model parameters and characterize the model uncertainty utilizing appropriate stochastic parameterizations.
\end{enumerate}

%Figure \ref{Fig:Overview} contains a schematic description of the auto-learning framework.
It is worthwhile highlighting several significant features of this framework. First, the framework differs fundamentally from nonlinear regression in determining model structure. Typical nonlinear regression retains many nonphysical terms with small coefficients to minimize error. Sparse regression can eliminate some terms \citep{santosa1986linear, tibshirani1996regression}, but lacks physical constraints. In contrast, the causal inference exploits the causal relationship, which reflects the underlying physics, to automatically determine the model structure that provides a physically explainable parsimonious model \citep{almomani2020entropic, fish2021entropic, kim2017causation, almomani2020erfit, elinger2021information, elinger2021causation}.

Second, the framework utilizes relationships and insights from the established models as a starting point rather than fitting models directly to short observational data. Candidate functional forms for the model structure are constrained by incorporating physical dependencies learned from the existing models. Unlike the pure data-driven method, the framework imposes known dynamical couplings from the known knowledge to inform model derivation. Physical insights from the established models guide the exploration of model structures and dependencies.

Third, analytic formulae are available to determine both the model structure via causal inference and the parameters using maximum likelihood \citep{chen2020learning}. In other words, the auto-learning framework developed here is wholly explicit and highly efficient. It differs from many other machine learning methods that contain black boxes and require expensive training. Therefore, the auto-learning framework is easy to use and generalizable to many scientific topics.

%\begin{figure}[tp]
%\hspace*{-0cm}\includegraphics[width=1.0\textwidth]{Figs/Fig1_3_new.pdf}
%\caption{The diagram of knowledge-informed learning algorithm.}\label{Fig:Overview}
%\end{figure}

\begin{figure}[ht]
    \hspace*{-0cm}\includegraphics[width=1.0\textwidth]{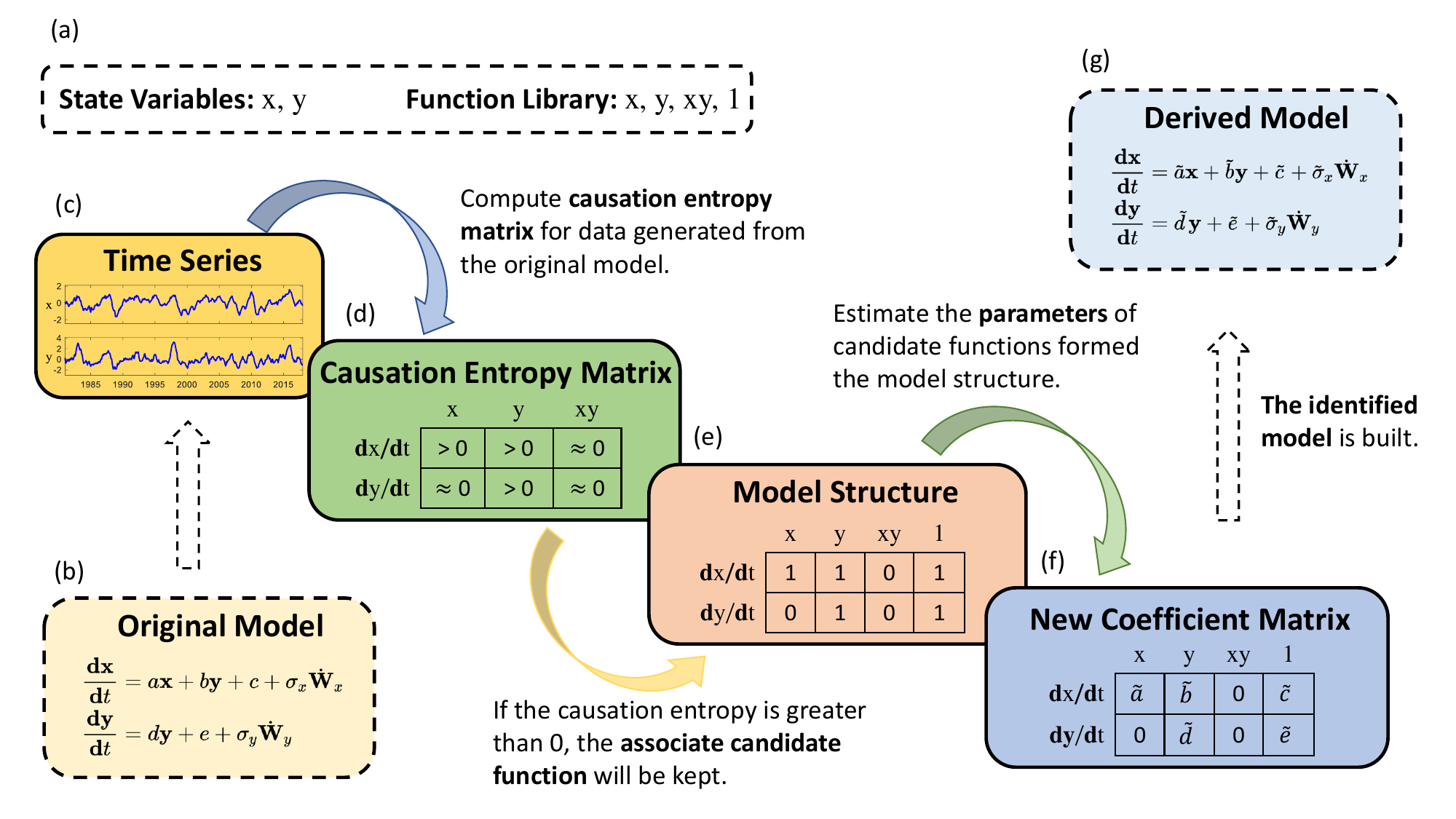}
    \caption{A schematic illustration of the auto-learning framework.   }\label{Fig:CEM}
\end{figure}

The details of Steps 2-5 are provided below, where a schematic illustration is shown in Figure \ref{Fig:CEM}. The established model used in Step 1, also known as the reference model, will be described at the end of this section.

\subsection{Determining model structure using causal inference}\label{Subsec:Causality}
\subsubsection{Determining the state variables}
The state variables of the model to be derived from the auto-learning framework are pre-determined. These variables are denoted by an $n$-dimensional column vector $\mathbf{U}=(u_1,\ldots, u_n)^{\mathtt{T}}$, where each $u_i$ is one state variable. To characterize ENSO diversity, $T_E$ and $T_C$ are the two variables appearing in all the models to be derived. Depending on the interests, additional variables are included in different models, such as the zonal current in the central Pacific to represent advection processes and the extra heat content variables to account for more ocean adjustment timescales.
\subsubsection{Developing a function library}
After determining the state variables of the target model, a library $\mathbf{f}$ consisting of in total $M$ possible candidate functions to describe the right-hand side of the target model is developed,
\begin{equation}\label{Library}
  \mathbf{f} = \{f_1,\ldots, f_{m-1}, f_m, f_{m+1}, \ldots, f_M\}.
\end{equation}
Typically, a large number of candidate functions is included in the library to allow for coverage over different possible dynamical features. Each $f_m$ is given by a linear or nonlinear function containing a few components of $\mathbf{U}$. In describing ENSO diversity, the library should at least include $T_E$, $T_C$, and some of their linear and nonlinear combinations, for example, $T_E^3$ and $T_C^3$. These functions may not appear in the resulting model from the auto-learning framework. But the library can provide as many potential candidate functions as possible that allow the causal inference to determine the useful ones.

\subsubsection{Determining model structure using causation inference}
Next, causal inference is utilized to evaluate the relevance of the candidate functions to the underlying dynamics and thus determine the model structure. To this end, the so-called causation entropy $C_{f_{m} \rightarrow \dot{u}_i \mid\left[\mathbf{f} \backslash {f}_{m}\right]}$ is computed to detect if the candidate function $f_m$ contributes to the right-hand side of the equation for each $u_i$, namely $\d u_i/\d t :=\dot{u}_i$. The causation entropy is given by \citep{elinger2021causation, almomani2020entropic, fish2021entropic}:
\begin{equation}\label{Causation_Entropy}
  C_{f_{m} \rightarrow \dot{u}_i \mid\left[\mathbf{f} \backslash {f}_{m}\right]}=H(\dot{u}_i|\left[\mathbf{f} \backslash {f}_{m}\right]) - H(\dot{u}_i|\mathbf{f}),
\end{equation}
where $\mathbf{f} \backslash {f}_{m}$ represents the set that contains all functions in $\mathbf{f}$ except $f_m$. In other words, $\mathbf{f} \backslash {f}_{m}$ contains $M-1$ candidate functions and is defined as
\begin{equation}\label{Library2}
  \mathbf{f}\backslash {f}_{m} = \{f_1,\ldots, f_{m-1},f_{m+1},\ldots, f_M\}.
\end{equation}
The term $H(\cdot|\cdot)$ is the conditional entropy, which is related to Shannon's entropy $H(\cdot)$ and the joint entropy $H(\cdot,
\cdot)$. For two multi-dimensional random variables $\mathbf{X}$ and $\mathbf{Y}$ (with the corresponding states being $\mathbf{x}$ and $\mathbf{y}$), they are defined as \citep{cover1999elements}:
\begin{equation}\label{Entropies}
\begin{split}
  H(\mathbf{X}) &= -\int_x p(\mathbf{x})\log(p(\mathbf{x}))\d \mathbf{x},\\
  H(\mathbf{Y}| \mathbf{X}) &= -\int_\mathbf{x}\int_\mathbf{y} p(\mathbf{x},\mathbf{y})\log(p(\mathbf{y}|\mathbf{x}))\d \mathbf{y}\d \mathbf{x},\\
  H(\mathbf{X},\mathbf{Y}) &= -\int_\mathbf{x}\int_\mathbf{y} p(\mathbf{x},\mathbf{y})\log(p(\mathbf{x},\mathbf{y}))\d \mathbf{y}\d \mathbf{x}.
\end{split}
\end{equation}
On the right-hand side of \eqref{Causation_Entropy}, the difference between the two conditional entropies indicates the information in $\dot{u}_i$ contributed by the specific function $f_m$ given the contributions from all the other functions. Thus, it tells if $f_m$ provides additional information to $\dot{u}_i$ conditioned on the other potential terms in the dynamics. It is worthwhile to highlight that the causation entropy in \eqref{Causation_Entropy} is fundamentally different from directly computing the correlation between $\dot{u}_i$ and $f_m$, as the causation entropy also considers the influence of the other library functions. If both $\dot{u}_i$ and $f_m$ are caused by a common factor $f_{m^\prime}$ in the function library, then $\dot{u}_i$ and $f_m$  can be highly correlated. Yet, in such a case, the causation entropy $C_{f_{m} \rightarrow \dot{u}_i \mid\left[\mathbf{f} \backslash {f}_{m}\right]}$ will be zero as $f_m$ is not the causation of $\dot{u}_i$.

The causation entropy is computed from each of the candidate functions in $\mathbf{f}$ to each $\dot{u}_i$. Thus, there are in total $NM$ causation entropies, which can be written as a $N\times M$ matrix, called the causation entropy matrix. Note that the dimension of $\mathbf{X}$ in \eqref{Entropies} is $M$ when it is applied to compute the difference between two terms on the right-hand side of the causation entropy in \eqref{Causation_Entropy}. This implies that the direct calculation of the entropies in \eqref{Entropies} involves a high-dimensional numerical integration, which is a well-known computationally challenging issue \citep{bellman1961dynamic}. To circumvent the direct numerical integration, the entropy calculation approximates all the joint and marginal distributions as Gaussians. In such a way, the causation entropy can be approximated by
\begin{equation}
\label{Entropy_Gaussians}
\begin{split}
C_{\mathbf{Z} \rightarrow \mathbf{X} | \mathbf{Y}} &=H(\mathbf{X} | \mathbf{Y})-H(\mathbf{X} | \mathbf{Y}, \mathbf{Z}) \\
& = H(\mathbf{X},\mathbf{Y}) - H(\mathbf{Y}) - H(\mathbf{X},\mathbf{Y},\mathbf{Z}) + H(\mathbf{Y},\mathbf{Z})\\
& \approx \frac{1}{2} \ln(\operatorname{det}(\mathbf{R}_{\mathbf{X}\mathbf{Y}}))-\frac{1}{2} \ln(\operatorname{det}(\mathbf{R}_{\mathbf{Y}})) - \frac{1}{2} \ln(\operatorname{det}(\mathbf{R}_{\mathbf{X}\mathbf{Y}\mathbf{Z}}))
 +\frac{1}{2} \ln(\operatorname{det}(\mathbf{R}_{\mathbf{Y}\mathbf{Z}})),
\end{split}
\end{equation}
where $\mathbf{R}_{\mathbf{XYZ}}$ denotes the covariance matrix of the state variables $(\mathbf{X},\mathbf{Y},\mathbf{Z})$ and similar for other covariances. The notations $\ln(\cdot)$ and $\operatorname{det}(\cdot)$ are the logarithm of a number and determinant of a matrix, respectively. Note that the second equality in \eqref{Entropy_Gaussians} is derived by using the chain rule of conditional entropy.

The simple and explicit expression in \eqref{Entropy_Gaussians} based on the Gaussian approximation can efficiently compute the causation entropy. It allows the computation of the causation entropy with a moderately large dimension, sufficient for deriving conceptual models. It is worth noting that the Gaussian approximation may lead to certain errors in computing the causation entropy if the true distribution is highly non-Gaussian. Nevertheless, the primary goal is not to obtain the exact value of the causation entropy. Instead, it suffices to detect if the causation entropy $C_{f_{m} \rightarrow \dot{u}_i \mid\left[\mathbf{f} \backslash {f}_{m}\right]}$ is nonzero (or practically above a small threshold value). In most applications, if a significant causal relationship is detected in the higher-order moments, it is very likely in the Gaussian approximation. This allows us to efficiently determine the sparse model structure, where the exact values of the nonzero coefficients on the right-hand side of the model will be calculated via a simple maximum likelihood estimation to be discussed in the following, which will be elaborated in the following. Note that the Gaussian approximation has been widely applied to compute various information measurements and leads to reasonably accurate results \citep{majda2018model, tippett2004measuring, kleeman2011information, branicki2012quantifying}.

With the $N\times M$ causation entropy matrix in hand, the next step is determining the model structure. This can be done by setting up a threshold value of the causation entropy and retaining only those candidate functions with the causation entropies exceeding the threshold. Causation entropy allows the resulting model to contain only functions that significantly contribute to the dynamics and facilitates a sparse model structure. Sparsity is crucial to discovering the correct underlying physics and prevents overfitting \citep{ying2019overview, brunton2016discovering}. It will also guarantee the robustness of the model in response to perturbations and allow the model to apply to certain extrapolation tests.

\subsubsection{Parameter estimation}\label{Subsec:ParameterEstimation}
The final step is to estimate the parameters in the resulting model. The model of $\mathbf{U}$ can be written in the vector form:
\begin{equation}\label{Identified_Model}
  \frac{\d \mathbf{U}}{\d t} = \boldsymbol\Phi(\mathbf{U})+\boldsymbol{\sigma} \dot{\mathbf{W}}(t),
\end{equation}
where $\boldsymbol\Phi(\mathbf{U})$ is an $n$-dimensional column vector with each entry representing the summation of the remaining nonlinear candidate functions from the causal inference for the corresponding component of $\mathbf{U}$. The matrix $\boldsymbol{\sigma}\in \mathbb{R}^{n\times d}$ is the noise amplitude, and $\dot{\mathbf{W}}(t) \in \mathbb{R}^{d\times 1}$ is a white noise. Denote by $\boldsymbol\Theta\in \mathbb{R}^s$ a column vector containing all the parameters in the model with $s\gg n$ in a typical situation. The first term on the right-hand side of \eqref{Identified_Model} can be written as
\begin{equation}\label{RHS_Model}
  \boldsymbol\Phi(\mathbf{U}) = \mathbf{M} \boldsymbol\Theta + \mathbf{Q},
\end{equation}
where $\mathbf{M}\in\mathbb{R}^{n\times s}$ is a matrix, each entry of which is a function of $\mathbf{U}$ while $\mathbf{Q}\in \mathbb{R}^n$ is a column vector that depends on $\mathbf{U}$. Only $\mathbf{M}$ is multiplied by the parameters $\boldsymbol\Theta$. The parameters $\boldsymbol\Theta$ can be easily determined using a maximum likelihood estimator. Meanwhile, the noise coefficients $\boldsymbol{\sigma}$ are computed. See \citep{chen2020learning} for the technical details. Notably, the entire parameter estimation can be solved via closed analytic formulae, making the procedure efficient and accurate. %Finally, energy-conserving constraints to the parameter values are often included in the parameter estimation procedure for quadratic nonlinear terms. These constraints in the conceptual model prevent the finite time blowup of the solution in the derived model and is physically consistent \citep{majda2012physics, harlim2014ensemble}. Remarkably, closed analytic formulae are still available for parameter estimation with such constraints.

\subsection{Auto-learning with latent variables}\label{Subsec:Latents}
Determining the state variables is a prerequisite for model identification. Yet, it is not always guaranteed that pre-determined state variables are appropriate for the entire system to capture the target features. Therefore, it is of practical importance to systematically recognize the sources of the model deficiency and suggest appropriate additional state variables that improve the model performance by augmenting the existing model with these additional quantities.

These can be achieved by incorporating latent variables into the auto-learning framework. Latent variables are additional processes whose observed time series are not directly available. They may act as surrogates for the missing physics due to the lack of a complete understanding of nature. They may also effectively characterize the state variables at the unresolved scales or play the role of additional stochastic parameterizations.

\begin{figure}[ht]
    \hspace*{-0cm}\includegraphics[width=1.0\textwidth]{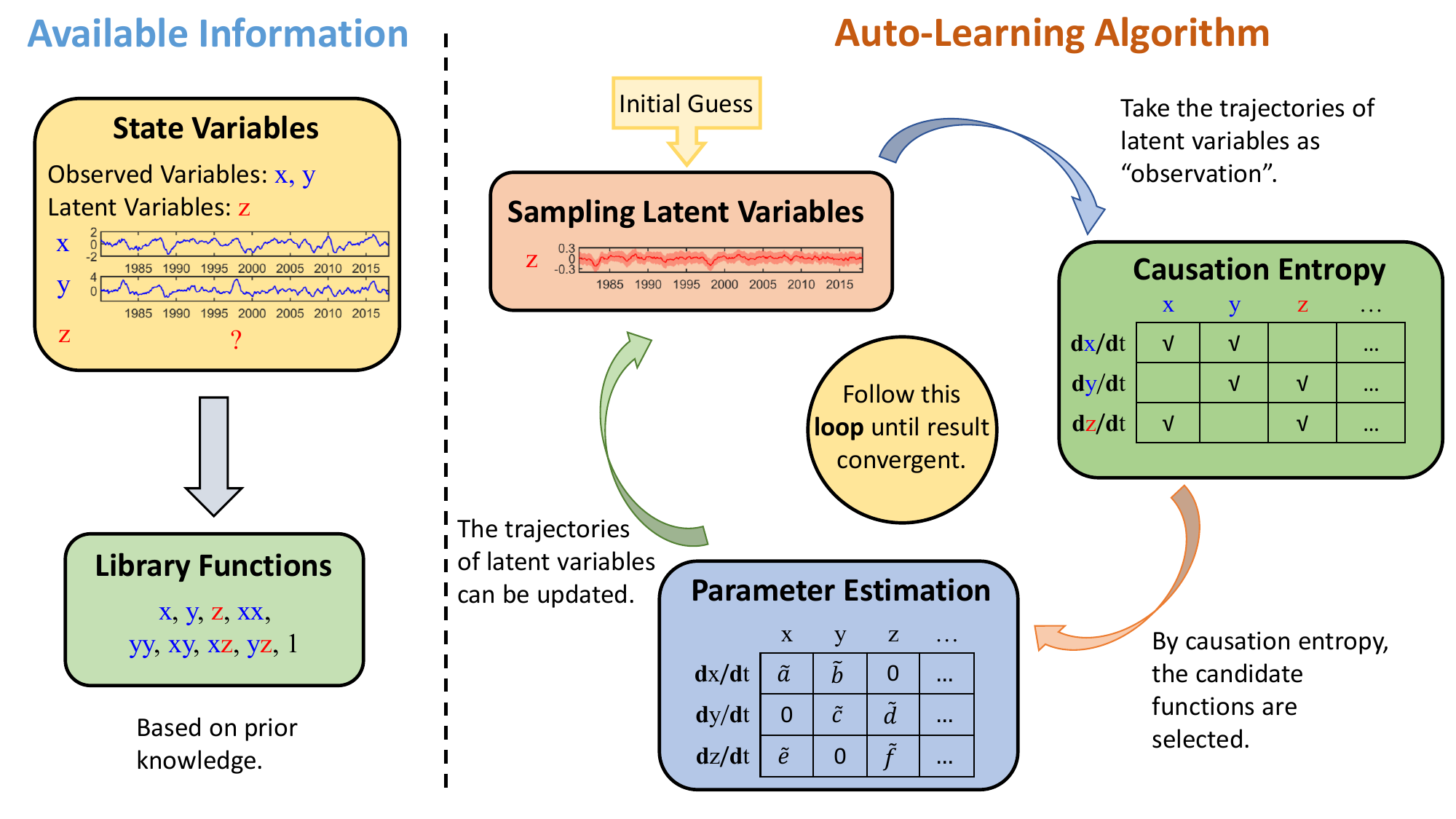}
    \caption{Overview of the auto-learning framework with latent variables.}
    \label{Fig:Overview_CELA}
\end{figure}

In the presence of latent variables, the auto-learning framework needs to be slightly modified by introducing a joint learning procedure to derive the equations for both groups of variables. Note that the time series is only available for the pre-determined state variables (also called observed variables in the following context). Therefore, an iterative method is used to alternate between inferring the time series of the latent variables and deriving the entire model structure. The model structure is determined using the causal inference discussed in Section \ref{Subsec:Causality} when the inferred time series of the latent variables is updated. The recovered time series of the latent variables, based on the current structure of the latent processes, is obtained by a Bayesian sampling technique \citep{chen2020learning}. The iterative procedure terminates when the model structure and parameters converge, which is ensured by the information revealed from observed variables. See Figure \ref{Fig:Overview_CELA} for the illustration of incorporating this iterative process into the auto-learning framework. The technical details of this iterative method can be found in \citep{chen2022causality}.

One crucial task is to explain the possible physical meanings of the latent variables once their time series is simulated or sampled from the augmented system. A stochastic conceptual model with a latent variable representing the variability in a fast time scale will be presented in Section \ref{Subsec:Latent}.

%The time series of the latent variables, Physical insights will be used to This will be presented in Section \ref{Subsec:Latent}, in which a stochastic conceptual model with latent variables will be derived. In such a case, the auto-learning framework needs to be slightly modified. In addition to the selected state variables, the model is augmented by a set of latent variables. Then the auto-learning framework involves a joint learning procedure to derive the equations for both groups of variables. Note that the time series is only available for the actual state variables (which are also called observed variables in the following context). Therefore, an iterative method is used to alternate between inferring the time series of the latent variables and deriving the model structure. The model structure is determined when the inferred time series of the latent variables is updated. The recovered time series of the latent variables, based on the current structure of the latent processes, is obtained by a Bayesian sampling technique \citep{chen2020learning}. The iterative procedure ends when the model structure and parameters converge, which is ensured by the information revealed from observed variables. See Figure \ref{Fig:Overview_CELA} for the illustration of incorporating this iterative process into the auto-learning framework. The technical details of this iterative method can be found in \citep{chen2022causality}.

\section{Reference Model and Observational data}\label{Sec:Observations}
\subsection{Observational data sets}
The following data sets are utilized in this study to provide some estimates of the observed ENSO properties, including its diversity. For the assessment of the model dynamical properties, monthly ocean temperatures, zonal current, thermocline depth and daily zonal winds at 850 hPa are used. For this, the oceanic variables are derived from the GODAS dataset \citep{behringer2004evaluation} (\url{https://www.ncdc.noaa.gov/oisst}), with thermocline depth calculated from potential temperature as the depth of the $20^o$C isotherm, while the atmospheric variable is obtained from the NCEP-NCAR reanalysis \citep{kalnay1996ncep} (\url{http://www.esrl.noaa.gov/psd/}). The analysis period is from 1982 to 2019. This period is used because it is during the satellite era and its observations are more reliable than those before. Anomalies are calculated by removing the monthly mean climatology of the whole period. For the assessment of the model performance on simulating different types of ENSO events, a longer period for SST (i.e., 1950-2020) is also used, which is obtained from the Extended Reconstructed Sea Surface Temperature version 5 \citep{huang2017extended} (\url{https://psl.noaa.gov/data/gridded/data.noaa.ersst.v5.html}). The Ni\~no 4 ($T_C$ in the model) and Ni\~no 3 ($T_E$ in the model) indices are the averages of SST anomalies over the regions 160$^o$E-150$^o$W, 5$^o$S-5$^o$N (CP) and 150$^o$-90$^o$W, 5$^o$S-5$^o$N (EP), respectively. The $h_W$ index in the model is the mean thermocline depth anomaly over the WP region (120$^o$E-180$^o$, 5$^o$S-5$^o$N). The $u$ index in the model is the mean mixed layer zonal current in the CP region.

\subsection{The reference model}\label{Sec:reference model}
Instead of using short observational data, an appropriate physics-informed reference model is adopted to generate a long time series that will be used as the input data for the auto-learning framework. This reference model captures many observed features and is calibrated using observational data. Therefore, it can be regarded as an appropriate surrogate for observations. Various new models will be derived using the auto-learning framework with simulated data from this model in Section \ref{Sec:Model}.

The reference model is given by a recently developed three-region multiscale stochastic model \citep{chen2022multiscale}. This model includes intraseasonal, interannual, and decadal processes and captures many observed dynamical and statistical features of ENSO diversity. The model is as follows,
\begin{subequations}\label{Reference_model}
    \begin{align}
    & \frac{\d u}{\d t}=-r u-\frac{a_1 b_0 \mu}{2}\left(T_C+T_E\right)+\beta_u \tau+\sigma_u \dot{W}_u, \label{R_u}\\
    & \frac{\d h_W}{\d t}=-r h_W-\frac{a_2 b_0 \mu}{2}\left(T_C+T_E\right)+\beta_h \tau+\sigma_h \dot{W}_h, \label{R_h}\\
    & \frac{\d T_C}{\d t}=\left(\frac{\gamma b_0 \mu}{2}-c_1\left(T_C,t\right)\right) T_C+\frac{\gamma b_0 \mu}{2} T_E+\gamma h_W+c_I Iu+C_u+\beta_C \tau+\sigma_C \dot{W}_C, \label{R_c}\\
    & \frac{\d T_E}{\d t}=\gamma h_W+\left(\frac{3 \gamma b_0 \mu}{2}-c_2(t)\right) T_E-\frac{\gamma b_0 \mu}{2} T_C+\beta_E \tau+\sigma_E \dot{W}_E, \label{R_e}\\
    & \frac{\d \tau}{\d t}=-d_\tau \tau+\sigma_\tau\left(T_C\right) \dot{W}_\tau, \label{R_tau}\\
    & \frac{\d I}{\d t}=-\lambda(I-m)+\sigma_I(I) \dot{W}_I \label{R_I}.
    \end{align}
\end{subequations}
Here, the interannual components, i.e., \eqref{R_u}--\eqref{R_e}, depict the main dynamics for ENSO. The intraseasonal variable $\tau$ in \eqref{R_tau} represents the amplitude of the random wind bursts \citep{thual2016simple} while the decadal variability $I$ in \eqref{R_I} represents the strength of the background Walker circulation. In the model, $T_C$ and $T_E$ are the SSTs in the CP and EP, $u$ is the zonal ocean current in the CP, and $h_W$ is the thermocline depth in the WP. As was discussed in \citep{chen2022multiscale}, $I$ also stands for the zonal SST difference between the WP and CP, which affects the interannual component by modulating the efficiency of the zonal advection through $Iu$ in \eqref{R_c}.

In this model, stochasticity plays a crucial role in coupling variables at different time scales and parameterizing the unresolved features in the model. First, intraseasonal wind bursts $\tau$ are modeled by a stochastic differential equation with state-dependent $\sigma_\tau(T_C)$. The wind burst $\tau$ is then coupled to the processes of the interannual part serving as an external forcing. In addition, four Gaussian random noises $\sigma_u\dot{W}_u$, $\sigma_h\dot{W}_h$, $\sigma_C\dot{W}_C$ and $\sigma_E\dot{W}_E$ are further added to the processes, which effectively parameterize the untracked factors, such as the subtropical atmospheric forcing at the Pacific Ocean. More broadly, the noise increases variability for statistics to match observations \citep{palmer2009stochastic}. Second, since the details of the background Walker circulation consist of uncertainties and randomness \citep{chen2015strong}, a simple but effective stochastic process is used to describe the temporal evolution of the decadal variability $I$ \citep{yang2021enso}. The multiplicative noise in the process of $I$ is aimed at guaranteeing its positivity because the long-term average of the background Walker circulation is non-negative. Besides, the effects of seasonality are added to both the wind activity and the collective damping to depict the seasonal phase-locking characteristics realistically, which manifests as the tendency of ENSO events to peak during boreal winter \citep{tziperman1997mechanisms, stein2014enso, fang2021effect}.

The model parameters are listed in Table \ref{tab:reference model}.

\begin{table}[ht]
    \centering
    \begin{tabular}{cccc}
    \hline \hline
    $[h]$ & $150 \mathrm{~m}$ & {$[T]$} & $7.5^{\circ} \mathrm{C}$ \\
    {$[u]$} & $1.5 \mathrm{~m} \mathrm{~s}^{-1}$ & {$[t]$} & 2 months \\
    {$[\tau]$} & $5 \mathrm{~m} \mathrm{~s}^{-1}$ & $d_\tau$ & 2 (1 month$^{-1}$) \\
    $\gamma$ & 0.45 & $r$ & 0.15 \\
    $\alpha_1$ & 0.0375 & $\alpha_2$ & 0.075 \\
    $b_0$ & 2.5 & $\mu$ & 0.5 \\
    $\sigma$ & $0.12 I$ & $\lambda$ & 2/60 (5 years$^{-1}$)\\
    $p(I)$ & 0.25 $\text{in}$ $I \in(0,4)$ & $c_U$ & 0.018 \\
    $\Phi(x)$ & $\int_b^x(y-m) p(y) d y$ & $\sigma_I(I)$ & $\sqrt{\frac{2}{p(I)}[-\lambda \Phi(I)]}$ \\
    $\beta_E$ & $0.1239\times(2-0.2 I)$ & $\beta_u$ & $-0.2 \beta_E$ \\
    $\beta_h$ & $-0.4 \beta_E$ & $\beta_C$ & $0.8 \beta_E$ \\
    $\sigma_u$ & 0.0310 & $\sigma_h$ & 0.0155 \\
    $\sigma_C$ & 0.0310 & $\sigma_E$ & 0.0232 \\
    $\sigma_\tau\left(T_C, t\right)$ & $0.9\left[\tanh \left(4.5 T_C\right)+1\right]\left[1+0.25 \cos \left(\frac{2 \pi}{6} t\right)\right]$ \\
    $c_1\left(T_C, t\right)$ & {$\left[15.6\left(T_C+\frac{0.75}{7.5}\right)^2+0.57\right]\left[1+0.4 \sin \left(\frac{2 \pi}{6} t\right)\right]$} \\
    $c_2(t)$ & {$0.9\left[1+0.4 \sin \left(\frac{2 \pi}{6} t+\frac{2 \pi}{6}\right)+0.2 \sin \left(\frac{2 \pi}{3} t+\frac{2 \pi}{6}\right)\right]$} \\
    \hline
    \end{tabular}
    \caption{Summary of the non-dimensional units and reference model parameters. The square brackets around the $h$, $u$, $\tau$, $T$, and $t$ mean the unit of each variable. Key parameters include the seasonally modulated damping rates $c_1$ and $c_2$, where $c_1$ has a $T_C$ nonlinearity justified by a heat budget analysis to capture asymmetric CP statistics. The wind burst noise $\sigma_{\tau}$ also depends on $T_C$ to generate irregularity and non-Gaussian SST distributions which is known as a multiplicative noise forcing and represents the effect of interannual SST anomalies on wind bursts \citep{jin2007ensemble}. $\sigma_I$ defines stochastic decadal variability for $I$, and the parameter $m$ is the mean of $I$ which can be computed from its PDF. }
    \label{tab:reference model}
\end{table}

\subsection{Definitions of different types of the ENSO events}\label{Sec:definitions}
The definitions of different El Ni\~no and La Ni\~na events are based on the average SST anomalies during boreal winter, i.e., from December to February (DJF). Following the definition in \citep{kug2009two}, EP events are defined as situations where the EP is warmer than the CP with EP SST anomalies above 0.5 $^o$C. Furthermore, based on the definitions in \citep{wang2019historical}, an extreme El Ni\~no event corresponds to a maximum EP SST anomaly above 2.5$^o$C between April and the following March. Similarly, CP events are defined as situations where the CP is warmer than the EP, and CP SST is above 0.5 $^o$C. A La Ni\~na event is also defined as either the CP or EP SST anomaly is below -0.5 $^o$C. Finally, when El Ni\~no or La Ni\~na occurs for two or more consecutive years, it is defined as a multi-year El Ni\~no or La Ni\~na, respectively.

\section{Validation Metrics for Assessing the Model Performance}\label{Sec:Metrics}
\begin{figure}[ht]
    \hspace*{-0cm}\includegraphics[width=1.0\textwidth]{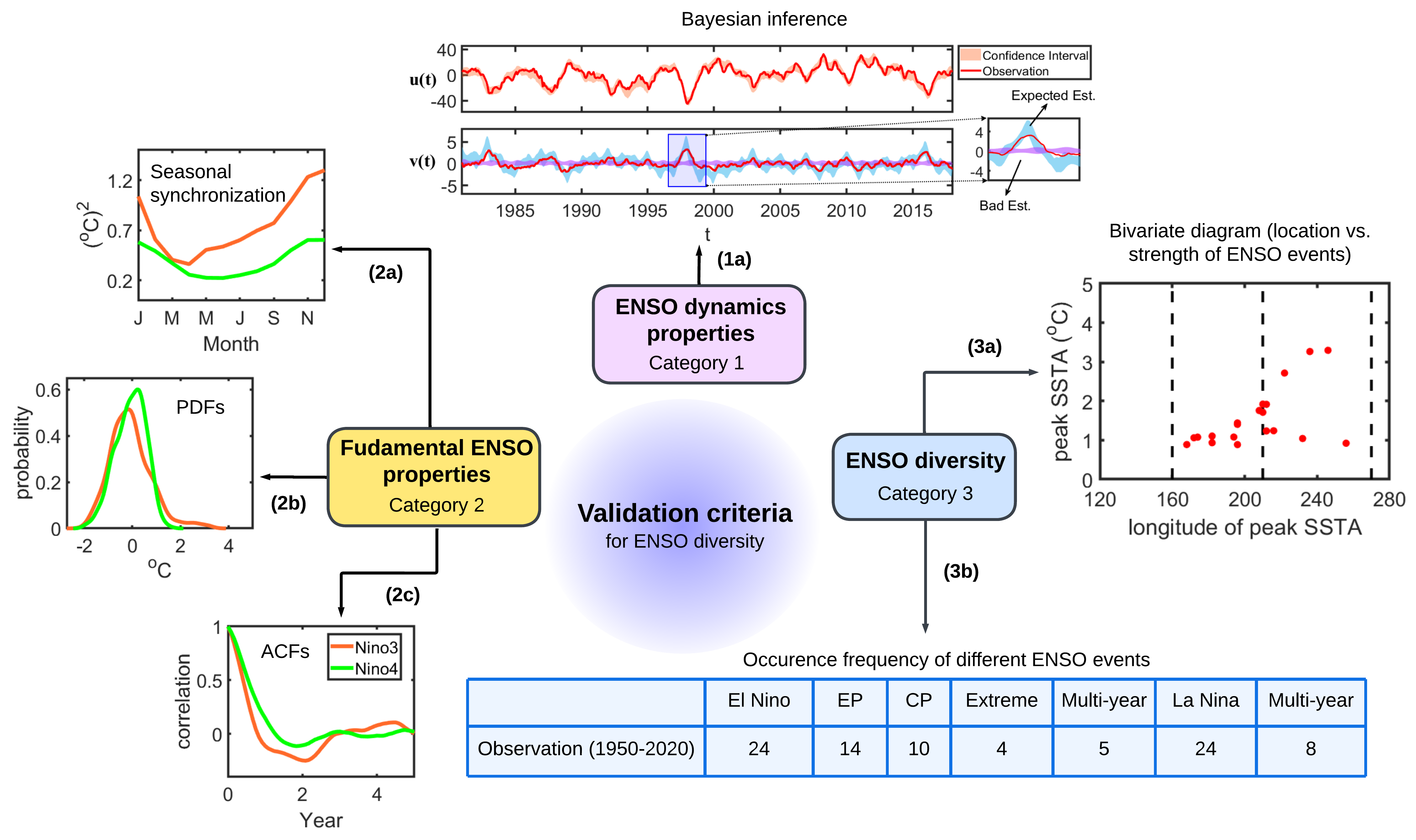}
    \caption{A summary of validation metrics for assessing the model performance.}\label{Fig:Category}
\end{figure}

The models derived from the auto-learning framework will be validated by a systematic set of metrics with observational data. These metrics measure the model performance in characterizing the essential features of ENSO diversity from both the dynamical and the statistical perspectives. The validation metrics are partitioned into three categories, which are summarized in Figure \ref{Fig:Category}. They will be discussed in detail below.

It is worth mentioning that the auto-learning framework does not explicitly exploit any of these metrics in deriving the stochastic conceptual models. Therefore, these metrics are objective for validating the performance of the resulting models.

\subsection{A metric for assessing the fundamental dynamical properties}
%\vspace{-0.1cm}
\subsubsection*{(1a) Nonlinear dependence between state variables}\label{Subsec:criterion_1a}

One of the main goals of developing a physical model is to reveal the underlying dynamics. This can, for instance, be achieved by quantifying the nonlinear inter-dependency between different state variables. Standard model evaluation metrics such as the lagged correlation are linear and may not properly reveal such nonlinear relationships.

Instead, a Bayesian inference approach is utilized. This idea is to test the ability of the model to reconstruct a subset of its state variables $\mathbf{u}(t)$ given the remaining subset $\mathbf{v}(t)$, using a Bayesian approach. Further denote the observational time series of these two sets of variables by $\mathbf{u}^{obs}(t)$ and $\mathbf{v}^{obs}(t)$. Loosely speaking, the Bayesian inference treats the time series of $\mathbf{v}^{obs}(t)$ as the constraint for the model to infer the time series $\mathbf{u}(t)$ and is then compared with $\mathbf{u}^{obs}(t)$. Technically, this can be achieved by computing the conditional distribution $p(\mathbf{u}(t)|\mathbf{v}^{obs}(t))$ at each time $t$ using the Bayes formula,
\begin{equation}\label{DA_formula}
  p(\mathbf{u}(t)|\mathbf{v}^{obs}(t)) \sim p(\mathbf{u}(t))\cdot p(\mathbf{v}^{obs}(t)|\mathbf{u}(t)),
\end{equation}
where $p(\mathbf{u}(t))$ is obtained from the model simulation, while $p(\mathbf{v}^{obs}(t)|\mathbf{u}(t))$ is computed based on both the model output and the available observation $\mathbf{v}^{obs}(t)$. In this work, a nonlinear Kalman smoother \citep{emerick2013ensemble, evensen2003ensemble} is utilized to recover $p(\mathbf{u}(t)|\mathbf{v}^{obs}(t))$. Then the mean values of $p(\mathbf{u}(t)|\mathbf{v}^{obs}(t))$ at different time instants provide the most likely time evolution of the model, given the observational time series. This allows point-wise comparisons between model simulations and observations. If recovering $\mathbf{u}(t)$ conditioned on $\mathbf{v}^{obs}(t)$ closely matches $\mathbf{u}^{obs}(t)$, it indicates that the model accurately captures the nonlinear dependence between these two sets of state variables. The performance of a model in capturing the fundamental dynamical properties is then assessed by analyzing whether individual simulation events match observations. In addition to comparing individual events, the correlation coefficient between $\mathbf{u}(t)$ and $\mathbf{u}^{obs}(t)$ can be adopted as a simple metric to quantify such dynamical consistency.

\subsection{Metrics for fundamental properties of the ENSO}

ENSO exhibits fundamental properties that conceptual models should aim to reproduce, including its amplitude, asymmetry, seasonal phase locking, and characteristic timescales. This category of validation metrics focuses on examining model skill in reproducing these fundamental ENSO properties.

To examine the performance of the models in reproducing fundamental ENSO properties, a time series of 20000 years from each model is generated. Omitting the first 1000 years as possibly the burning period, the remaining time series of 19000 years is divided into 271 non-overlapping segments, each having a length of 70 years as the observational time series. These different segments allow us to estimate the 95\% confidence interval when computing these statistics.

\subsubsection*{(2a) Seasonal synchronization}

The seasonal variation of the SST variance is a standard statistic reflecting the seasonal synchronization of a model. This statistic is computed in both the EP and CP regions.

\subsubsection*{(2b) Probability density functions (PDFs)}
ENSO is asymmetrical. Its PDFs are positively skewed in the EP and negatively skewed in the CP. These asymmetry properties are caused by ENSO extreme events. These extreme events, known as the super El Ni\~no, have become more frequent in recent decades and have significantly impacted the global climate \citep{thual2019statistical, chen2017formation, hameed2018model}. The super El Ni\~no usually occurs in the EP region while most CP events are moderate with almost no extreme CP El Ni\~no. Therefore, the SST PDFs of these two regions are crucial for examining the performance of the model in capturing ENSO extreme events.

\subsubsection*{(2c) Autocorrelation functions (ACFs)}
The ACF measures the overall memory of a chaotic system and describes the averaged convergence rate of the statistics toward the climatology \citep{gubner2006probability}. It is a valuable quantity to characterize the transient behavior and directly affects the model predictability. If the ACF of a model is strongly biased, then the model will fail to describe and predict any ENSO events, even if it has perfectly non-Gaussian PDFs.

Note that despite being widely used for many other studies of time series analysis, the power spectrum is not included as a validation metric in this paper since the information delivered by the spectrum overlaps with that of the temporal correlations \citep{priestley1981spectral, kay1988modern,percival1993spectral}.

\subsection{Metrics for ENSO diversity}
Besides the fundamental properties of the ENSO, two additional quantitative metrics are designed here to evaluate model fidelity in reproducing the observed spatial pattern diversity of ENSO categories and proportions.

\subsubsection*{(3a) Bivariate diagram of location and strength of ENSO events}
The first criterion in this set is a bivariate diagram that depicts the location and the strength of the ENSO events, which was first proposed in \citep{capotondi2015understanding}. This is an essential criterion to characterize the spatiotemporal pattern of the ENSO events.

This metric requires the spatiotemporal reconstruction of the SST data. Denote by $x$ and $t$ the longitude and the time, $\mbox{SST}(x,t)$ can be obtained from a multivariate linear regression based on the two time series of the CP SST $T_C$ and the EP SST $T_E$ from the model:
\begin{equation}\label{Regression}
  \mbox{SST}(x,t) = r_C(x)T_C(t) + r_E(x)T_E(t),
\end{equation}
where $r_C$ and $r_E$ are the two regression coefficients between the observational time series at the longitude $x$ and the $T_C$ and $T_E$ time series, respectively.

\subsubsection*{(3b) Occurrence frequency of different types of events}
The second criterion is the occurrence frequency of different ENSO events. This is a refined metric that counts the number of events in multiple sub-categories, including all the El Ni\~no events, the EP El Ni\~nos, the CP El Ni\~nos, the extreme El Ni\~nos, the multi-year El Ni\~nos, all the La Ni\~na events, and the multi-year La Ni\~nas \citep{chen2022multiscale} according to the definitions mentioned in Section \ref{Sec:definitions}.

\section{A Hierarchy of Stochastic Conceptual Models for ENSO Diversity}\label{Sec:Model}

\subsection{Overview of the models}
The auto-learning framework is applied to derive five models with different state variables and degrees of freedom, summarized in Table \ref{Table:Overview_Models}. The detailed expressions of these models are included in the Appendix. Depending on the prescribed state variables in the model, the corresponding time series generated from the reference model is utilized for model identification in the auto-learning framework.

\begin{table}[ht]\begin{center}
\begin{tabular}{|l|l|l|l|}
  \hline
  % after \\: \hline or \cline{col1-col2} \cline{col3-col4} ...
  Model Name  & Acronym & Dimension & State Variables \\\hline
  Interannual Intraseasonal Decadal Model with Advection & $I_a I_s D M A$ & 6 & $u, h_W, T_C, T_E, \tau, I$ \\
  Interannual Intraseasonal Model with Advection & $I_a I_s M A$ & 5 & $u, h_W, T_C, T_E, \tau$ \\
  Interannual Intraseasonal Decadal Model & $I_a I_s D M$ & 5 & $h_W, T_C, T_E,\tau, I$ \\
  Interannual Intraseasonal Model & $I_a I_s M$ & 4 & $h_W, T_C, T_E, \tau$ \\
  Interannual Model & $I_a M$ & 3 & $h_W, T_C, T_E$ \\
  \hline
\end{tabular}\end{center}\caption{Summary of the five models derived from the auto-learning framework}\label{Table:Overview_Models}
\end{table}

\begin{enumerate}
    \item \textbf{Interannual Intraseasonal Decadal Model with Advection ($I_a I_s D M A$)}: This 6-dimensional model involves three time scales, consistent with the multiscale system in \citep{chen2022multiscale}. 
    %According to Part I of this work, it has a slightly more parsimonious structure than the reference model and can capture most of the desirable features of ENSO diversity features.
    \item \textbf{Interannual Intraseasonal Model with Advection ($I_a I_s M A$)}: This 5-dimensional model excludes the decadal variable $I$ and only contains interactions between the interannual and intraseasonal variables. This model will be used to examine if the decadal modulation significantly impacts the dynamics and statistics of ENSO diversity.
    \item \textbf{Interannual Intraseasonal Decadal Model ($I_a I_s D M$)}: This 5-dimensional model involves three time scales but excludes the variable describing the ocean zonal current. This current variable was introduced to resolve advection, which plays a strong role in the central Pacific \citep{vialard2001model}. However, \citep{jin1999thermocline} pointed out that advection is implicitly included in the standard Recharge Oscillator framework, suggesting that the $u$ variable may not be necessary.
    \item \textbf{Interannual Intraseasonal Model ($I_a I_s M$)}: This 4-dimensional model excludes the decadal variability and the ocean zonal current. It aims to use a minimum number of state variables to characterize the interactions between interannual and intraseasonal variabilities.
    \item \textbf{Interannual Model ($I_a M$)}: This 3-dimensional model further excludes the intraseasonal wind variable. Therefore, it reduces to a single time scale model describing the interannual variabilities. Note that this model still includes stochastic forcing. However, the intraseasonal variability, which is characterized by a state-dependent red noise is omitted in this model. This model will be compared with the $I_a I_s M$ to reveal the crucial role of intraseasonal variability in contributing to ENSO diversity \citep{chen2015strong, eisenman2005westerly, jin2007ensemble}. Despite the lack of skill to fully characterize ENSO diversity, this $I_a M$ will serve as the basis for the auto-learning framework to include latent variables for model improvement.
\end{enumerate}

The entire library of candidate functions includes all the linear and quadratic functions of the six state variables, the cubic function of each variable, and a constant term representing the forcing:
\begin{equation}\label{6d_candidate}
    \begin{gathered}
    u,\ h_W,\ T_C,\ T_E,\ \tau,\ I,\ u^2,\ h_W^2,\ T_C^2,\ T_E^2,\ \tau^2,\ I^2, \\
    uh_W,\ uT_C,\ uT_E,\ u\tau,\ uI,\ h_WT_C,\ h_WT_E,\ h_W\tau,\ h_WI, \\
    T_CT_E,\ T_C\tau,\ T_CI,\ T_E \tau,\ T_E I,\ \tau I,\ u^3,\ h_W^3,\ T_C^3,\  T_E^3,\ \tau^3,\ I^3,\ 1.
    \end{gathered}
\end{equation}
For models with only a subset of the six possible variables, the corresponding candidate functions involving the missing state variables are removed from the library. %Note that the functions of multiplicative noises are not included in the library since they could not be identified from the framework. The multiplicative noise of $\tau$ and $I$ are assumed to be known and they are the same as the reference model.

It is worth mentioning that the derivation of the 6-dimensional interannual intraseasonal decadal model with advection ($I_a I_s D M A$) is also considered a proof-of-concept twin experiment to validate the auto-learning framework. As seen in the following (e.g., Figure \ref{Fig:Model_Structures_II}), the identified model has almost an identical structure to the reference model. The simulations from these two models also resemble each other. This evidence justifies the auto-learning framework.

\subsection{Intercomparison between different models}
To provide an overview of the model performance in characterizing different features of ENSO diversity, let us start with an intercomparison between these models.

\subsubsection{Criterion category 1: assessment of the fundamental dynamical properties}

\begin{figure}[H]
    \centering
    \includegraphics[width=1.0\textwidth]{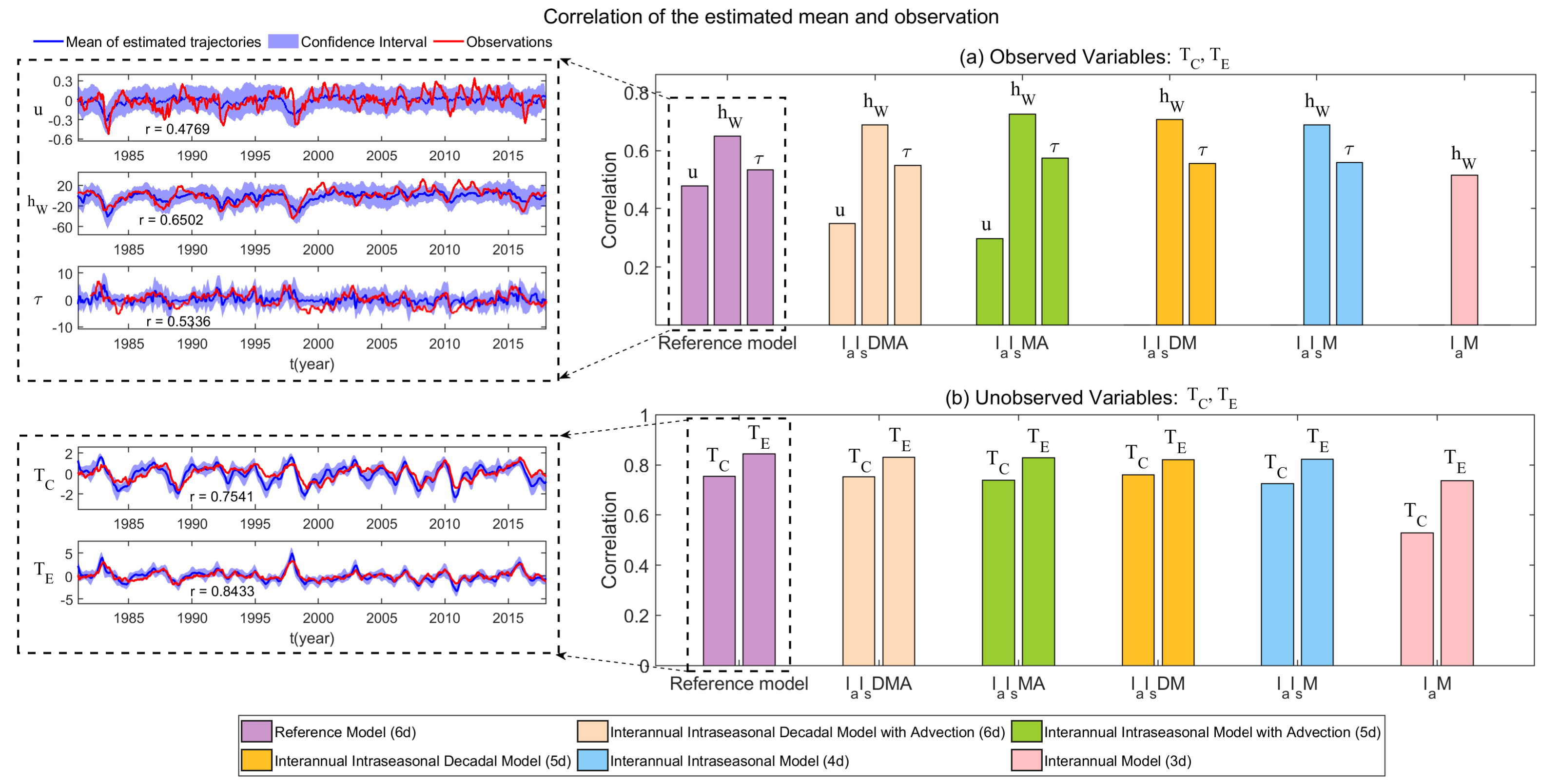}
    \caption{Comparison between different models using Criterion (1a) -- assessing the fundamental dynamical properties via Bayesian inference. The correlations between the recovered mean time series and observations are used to quantify the performance of these models in the right box chart. Panel (a): correlations of the results for recovery of $u$, $h_W$ and $\tau$ by observing $T_C$, $T_E$ and $I$. Panel (b): correlations of the results for recovery of $T_C$ and $T_E$ by observing other variables. On the left side, there is an example showing the recovered time series of the reference model in two experiments. The red lines represent the observation, the solid blue lines represent the ensemble mean simulation, and the blue shading areas are the confidence intervals of the ensemble simulations. The correlation $r$ between the truth (red line) and the recovered ensemble mean time series (solid blue line) for each state variable is also listed. }
    \label{fig:PartII_DA_new}
\end{figure}

Figure \ref{fig:PartII_DA_new} shows the results using the Bayesian inference to access the fundamental dynamical properties. The correlation between the recovered mean time series and observations is used to quantify the performance of different models. Panel (a) and (b) illustrate the results of two experiments. In the first experiment, the SST variables $T_C$, $T_E$, and decadal variability $I$ are observed to recover the time series of the remaining variables in each model. Note that $I$ is the slowest variable and is always treated as observations. The results show that all models, except the 3-dimensional interannual model, provide comparably high correlations in recovering $h_W$, which are around 0.7. In contrast, the 3-dimensional interannual model only has a correlation around 0.5. One interesting finding in Panel (a) is that recovery of $u$ seems quite challenging because of the weak coupling relationship between the SST and current $u$ in the models. This also indicates that $u$ may only have a weak direct contribution to modulating the SST variables. It is consistent with the fact that $u$ is heavily driven by the atmospheric wind, and its role has primarily been included by the wind. Such a finding suggests that the process of $u$ may not be entirely necessary for building the minimum model of ENSO diversity.

Panel (b) shows the results of the second experiment. Here, the two SST variables $T_C$ and $T_E$ are unobserved, while the time series of all the remaining variables in each model are utilized as observations. All first four models (from the 6-dimensional interannual intraseasonal decadal model with advection to the 4-dimensional interannual intraseasonal model) perform comparably well in recovering the two SST variables, where the correlation is around 0.75 for $T_C$ and 0.82 for $T_E$, respectively. These values are similar to those of the reference model in the inference. In contrast, the correlations using the 3-dimensional interannual model reduce to 0.53 for $T_C$ and 0.73 for $T_E$. These numbers are still above the standard threshold for skillful results, with a correlation of 0.5. However, compared with other models, the 3-dimensional interannual is less accurate in characterizing the interdependence between different state variables.

Two dashed boxes on the left side in Figure \ref{fig:PartII_DA_new} show the direct results from Bayesian inference, including the ensemble mean time series (blue), the confidence interval, and the observations (red) of the reference model as an example. The reference model has a good performance in characterizing the nonlinear relationship between the state variables. Specifically, the recovered $T_E$ time series shows significant consistency with the observation in depicting the events around 1982-1983 and 1997-1998, and the 2015-2016 super El Ni\~nos.

%\clearpage

\subsubsection{Criterion category 2: quantification of fundamental ENSO properties}

\begin{figure}[ht]
    \hspace*{-0cm}\includegraphics[width=1.0\textwidth]{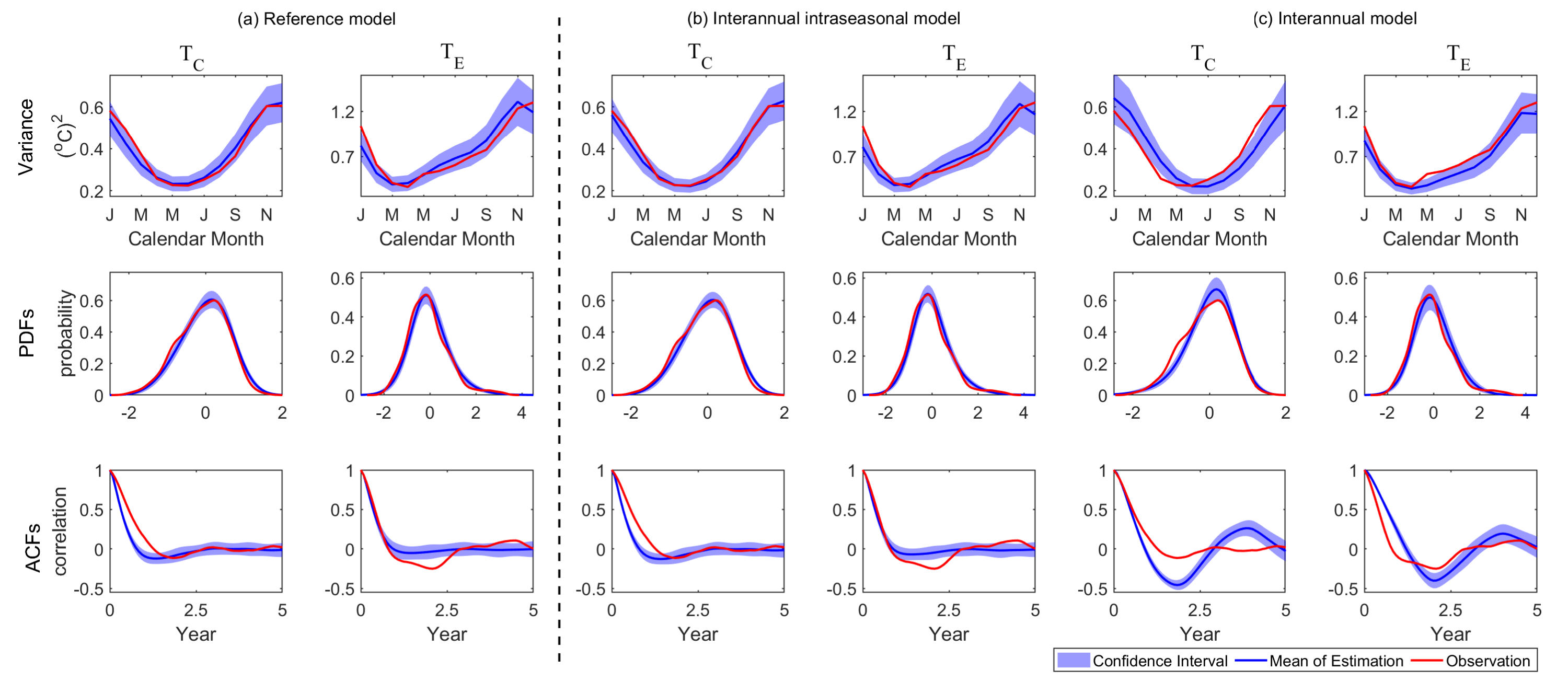}
    \caption{Comparison between different models using Criteria 2a--2c including the seasonal variation of the variance, PDFs, and ACFs of $T_C$ and $T_E$.}
    \label{Fig:PartII_Category2}
\end{figure}

The first row in Figure \ref{Fig:PartII_Category2} shows the seasonal dependence of the $T_C$ and $T_E$ amplitude (Criterion 2a). Only the results of three models are shown to avoid redundancy -- the first four models (from the 6-dimensional interannual intraseasonal decadal model with advection to the 4-dimensional interannual intraseasonal model) perform comparably similarly. All the models have the skill to reproduce seasonal synchronization. %It is worth highlighting that the observed SST variation is not a purely sinusoidal function over a year. Nevertheless, exploiting simple sinusoidal functions to parameterize the time-dependence of certain model parameters is sufficient for these models to accurately capture the observed seasonal synchronization since nonlinearity interacts with these sinusoidal functions.

\begin{figure}[ht]
    \hspace*{-0cm}\includegraphics[width=1.0\textwidth]{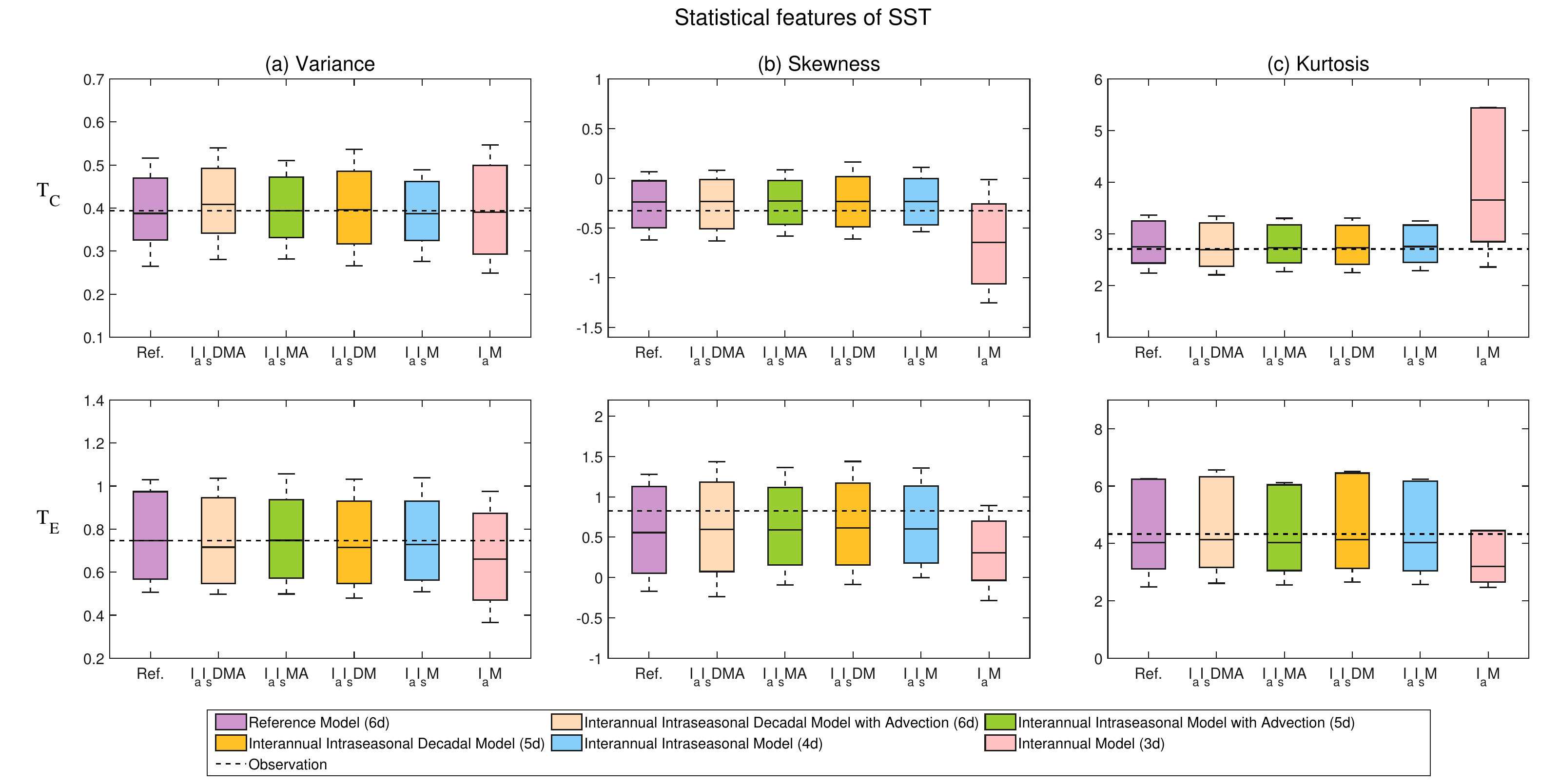}
    \caption{Additional supporting information of Criterion 2b. The second, third, and fourth moments, i.e., the variance, skewness, and kurtosis, of $T_C$ and $T_E$ for each model. These statistics can be computed using the PDFs in Figure \ref{Fig:PartII_Category2}. In each box plot, the central mark on each box denotes the median, while the bottom and top edges of the box signify the 5th and 95th percentiles, respectively. The horizontal dashed line in each panel shows the true value using observations. }
    \label{Fig:statistics_II}
\end{figure}

Next, the comparison of the non-Gaussian PDFs  (Criterion 2b) is included in the second row of Figure \ref{Fig:PartII_Category2}. All the models have similar profiles regarding these PDFs, especially the sign of the skewness for $T_C$ and $T_E$. The only noticeable difference among different models is that the PDF of $T_C$ associated with the 3-dimensional interannual model is slightly biased compared with other models. Yet, it is worth remarking that the higher-order moments that characterize non-Gaussian features may not be easily seen from the standard PDF plots. Therefore, Figure \ref{Fig:statistics_II} shows the box plots for the variance, skewness, and kurtosis from the model based on the 271 simulations. Despite the correct signs of the skewness being correctly recovered by all these models, the 3-dimensional interannual model has significant biases in reproducing the skewness and kurtosis values. For $T_C$, both the skewness and kurtosis are stronger than the observations. In particular, the kurtosis exceeds $3$, which means the probability tail is fatter than a Gaussian distribution that will cause more occurrence of extreme CP La Ni\~na events in such a model. In contrast, both the skewness and kurtosis are weaker than observations for $T_E$. This corresponds to the underestimation of the extreme EP El Ni\~no events. These findings indicate that the state-dependent red noise is crucial in reproducing these strong non-Gaussian features. The simple additive white noise is insufficient to capture the observed features in nature. 

Finally, the third row in Figure \ref{Fig:PartII_Category2} compares different models in recovering the ACFs (Criterion 2c). Similar to the variances and PDFs, only the 3-dimensional interannual model demonstrates different features compared with other models. The ACFs of the 3-dimensional interannual model have stronger oscillation patterns for both $T_C$ and $T_E$. Such a feature is related to the strong peak in the power spectrums, as the noise in the temporal direction is weaker than in the other models.

\subsubsection{Criterion category 3: evaluation of ENSO diversity}
%The last category of validation metrics aims to evaluate the statistical features of the dynamics.

\begin{figure}[ht]
    \hspace*{-0cm}\includegraphics[width=1.0\textwidth]{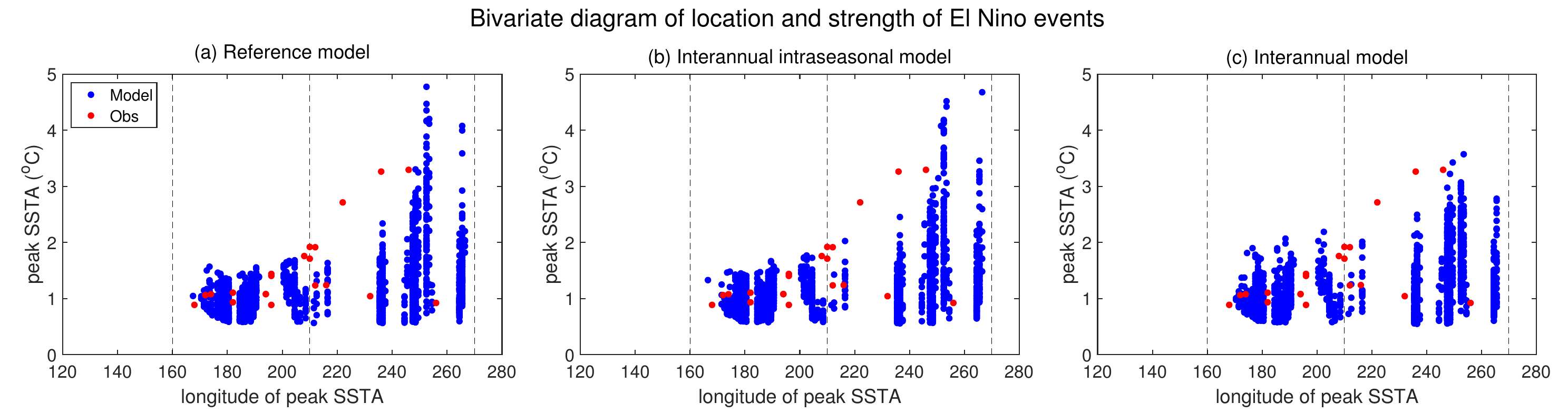}
    \caption{Comparison between different models using Criterion 3a: bivariate diagram of location and strength of El Ni\~no events. Each panel shows the distribution of equatorial Pacific SSTa maxima for El Ni\~no events of the standard simulations. For each of the El Ni\~no events, the winter-mean SST anomalies are averaged over the equatorial zone, and then the Pacific zonal maximum is located. For the presentation purpose, the blue dots show the El Ni\~no events based on the time series of only the first 4000 years of simulation in each panel. }
    \label{Fig:BPs_II}
\end{figure}

Figure \ref{Fig:BPs_II} shows the bivariate diagram of the location and strength of El Ni\~no events (Criterion 3a). Again, the first four derived models have similar performance. Thus, only the 4-dimensional interannual intraseasonal model and 3-dimensional interannual model are shown. A continuous spectrum of events across the equatorial Pacific is obtained for each model. In general, the CP events always exhibit weak to moderate amplitudes, while amounts of the strong events can be found in the EP region, which is consistent with observations \citep{capotondi2015understanding}. Among different models, the 3-dimensional interannual model is the only one that underestimates the number of extreme events in EP. This is related to the inaccuracy in the skewness and kurtosis shown in Figure \ref{Fig:statistics_II}.

\begin{figure}[ht]
    \hspace*{-0cm}\includegraphics[width=1.0\textwidth]{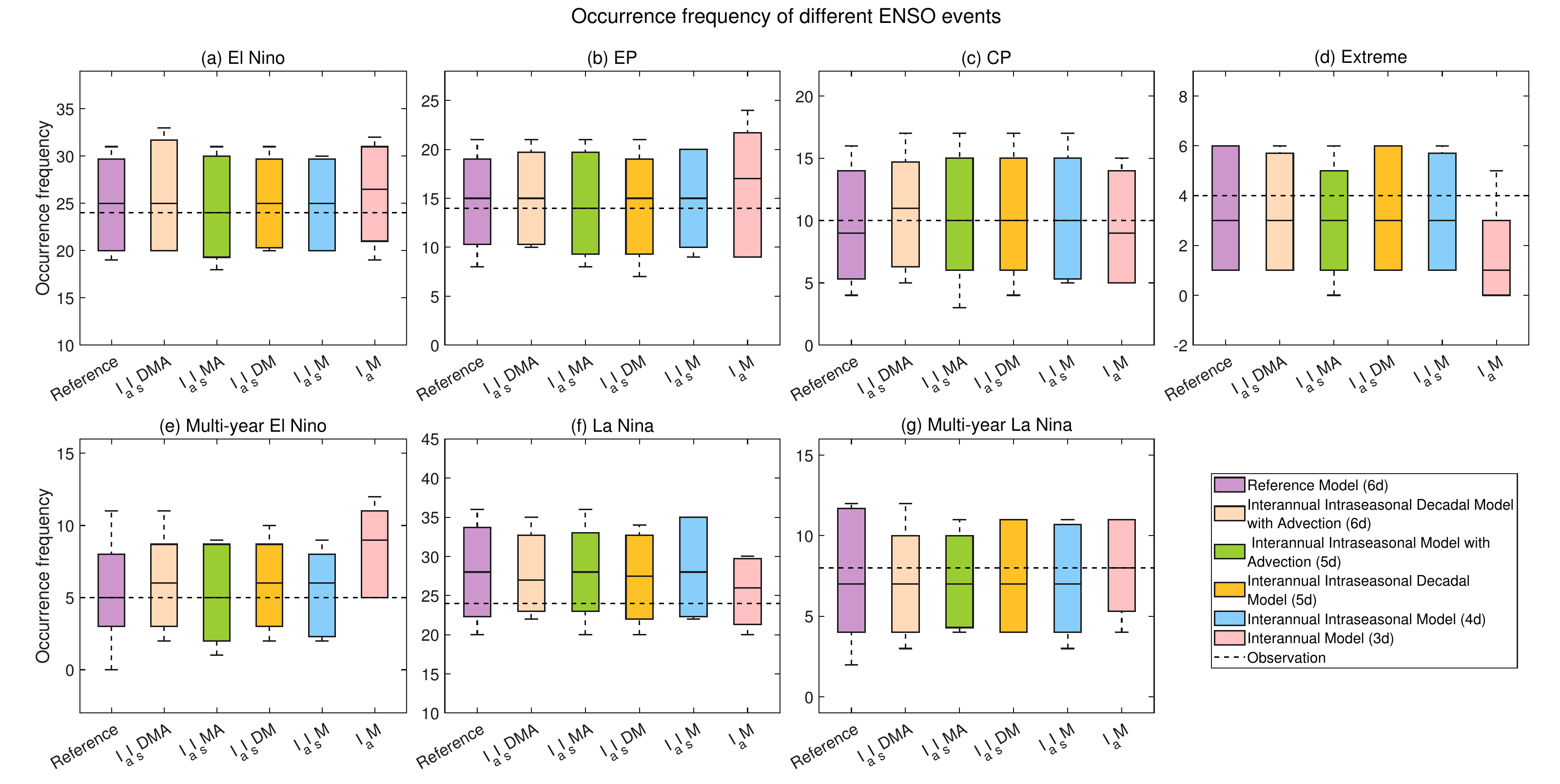}
    \caption{Comparison between different models using Criterion 3b: occurrence frequency of different types of events.  In the box plot, the central mark on each box denotes the median, while the bottom and top edges of the box signify the 5th and 95th percentiles, respectively. The horizontal dashed line in each panel shows the true value using observations. }
    \label{Fig:EC_II}
\end{figure}

Figure \ref{Fig:EC_II} presents the occurrence frequency of different types of events (Criterion 3b). All the models, except the 3-dimensional interannual model, can capture the observed occurrence frequency of all the seven types of ENSO events. This guarantees the skill of these models in simulating ENSO diversity and its complex features. Note that the 3-dimensional interannual model nevertheless reproduces the occurrence frequency of most ENSO events. The only two exceptions are the underestimation of extreme events, which has already been discussed in Figures \ref{Fig:statistics_II} and \ref{Fig:BPs_II}, and the overestimation of the number of multi-year El Ni\~no events due to the longer memory shown in the third row of Figure \ref{Fig:PartII_Category2}.

\clearpage

\section{Mechanisms and the Minimum Model for ENSO Diversity}\label{Sec:Discussion}

\begin{figure}[ht]
    \hspace*{-0cm}\includegraphics[width=1.0\textwidth]{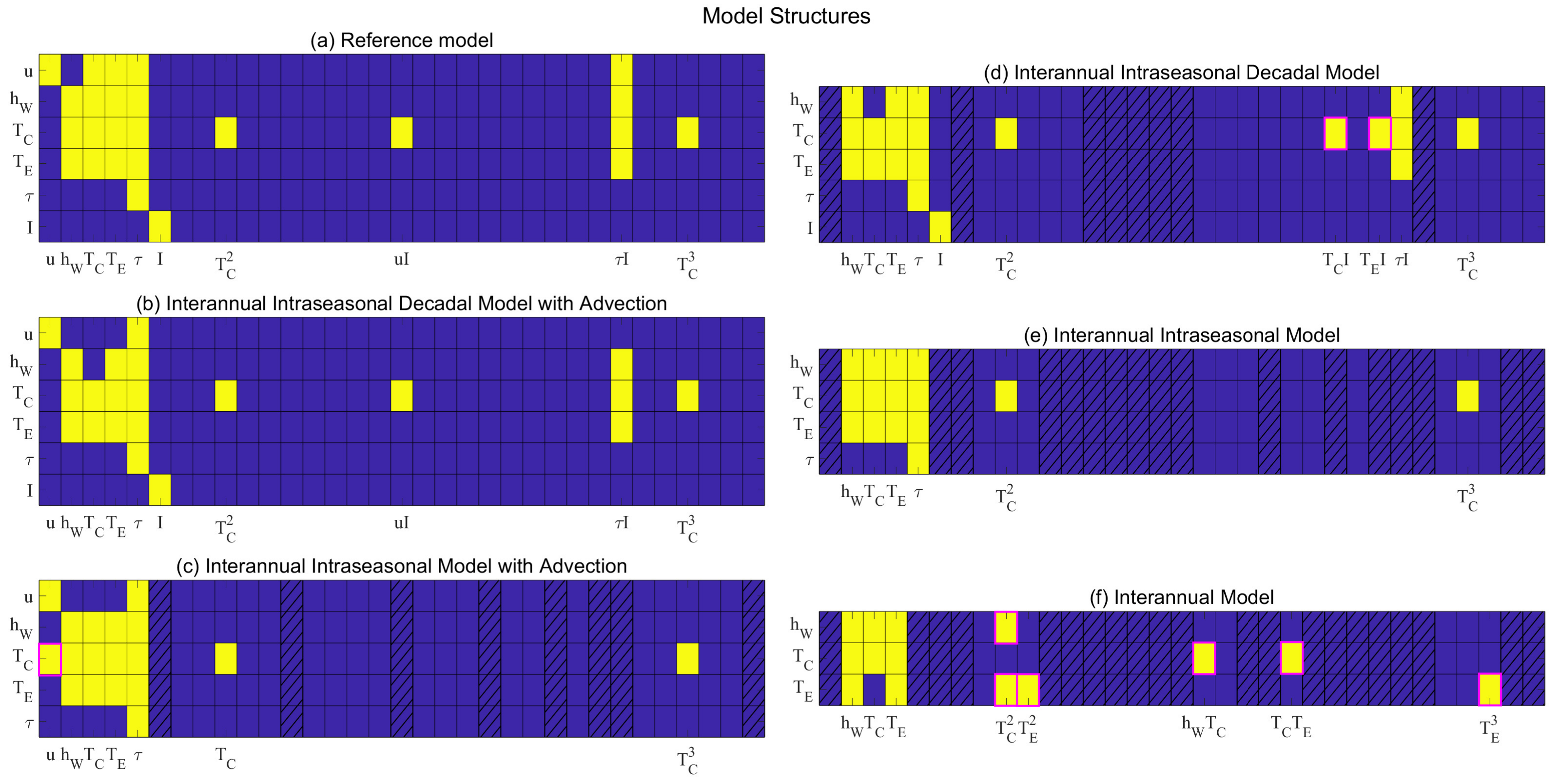}
    \caption{A summary of the model structures for the reference model (Panel (a)) and the five derived models from the auto-learning framework (Panels (b)--(f)). The yellow grids indicate the terms that appear in the model, while the blue ones mean the functions do not exist on the right-hand side of the derived model. The columns with dark blue shadings indicate the functions not built in the library for the corresponding model. The pink boxes mark the terms that are present in the derived models but not in the reference model.}
    \label{Fig:Model_Structures_II}
\end{figure}

Figure \ref{Fig:Model_Structures_II} illustrates the main structure of each model determined by the auto-learning framework. Note that the same number of columns is utilized for all the models to keep a unified format in the figure. The columns with dark blue shadings indicate the functions not built in the library for the corresponding model. Additional additive noise is added to the equations of $u$, $h_W$, $T_E$, and $T_C$, while multiplicative noise is imposed on the $\tau$ and $I$. The details of the models can be found in the Appendix.

\subsection{The minimum model}
The 4-dimensional $I_sI_aM$ is the minimum model that succeeds in capturing the target features of ENSO diversity characterized by the validation criteria in Figures \ref{fig:PartII_DA_new}--\ref{Fig:EC_II}. A further reduction of the model dimension to $3$ by eliminating the intraseasonal variability will significantly decrease the model's capacity to reproduce many observed features.

%The model, according to Panel (e) of Figure \ref{Fig:Model_Structures_II} or \eqref{IaIsM}, also has a relatively simple structure. It is almost a linear model containing the interactions between the three interannual variabilities. The only nonlinear terms are the cubic damping in the process of $T_C$ and a quadratic term $T_E\tau$ on the right-hand side of the $T_E$ process. The former has been shown to be crucial to reproduce the non-Gaussian features in $T_C$ (see Section \ref{Subsec:Linear}). The latter can be regarded as a multiplicative noise, where the wind forcing coefficient depends on $T_E$.

The model, according to Panel (e) of Figure \ref{Fig:Model_Structures_II} or \eqref{IaIsM}, also has a relatively simple structure. It is almost a linear model containing the interactions between the three interannual variabilities. The only nonlinear terms are the cubic damping in the process of $T_C$, which have been shown to be crucial in reproducing the non-Gaussian features in $T_C$ \citep{chen2022multiscale}. Similarly, the multiplicative noise in the wind burst process contributes to the extreme events and non-Gaussian features in $T_E$ \citep{thual2016simple, levine2017simple}. %(see Section \ref{Subsec:Linear}).

\begin{figure}[ht]
    \hspace*{-0cm}\includegraphics[width=1.0\textwidth]{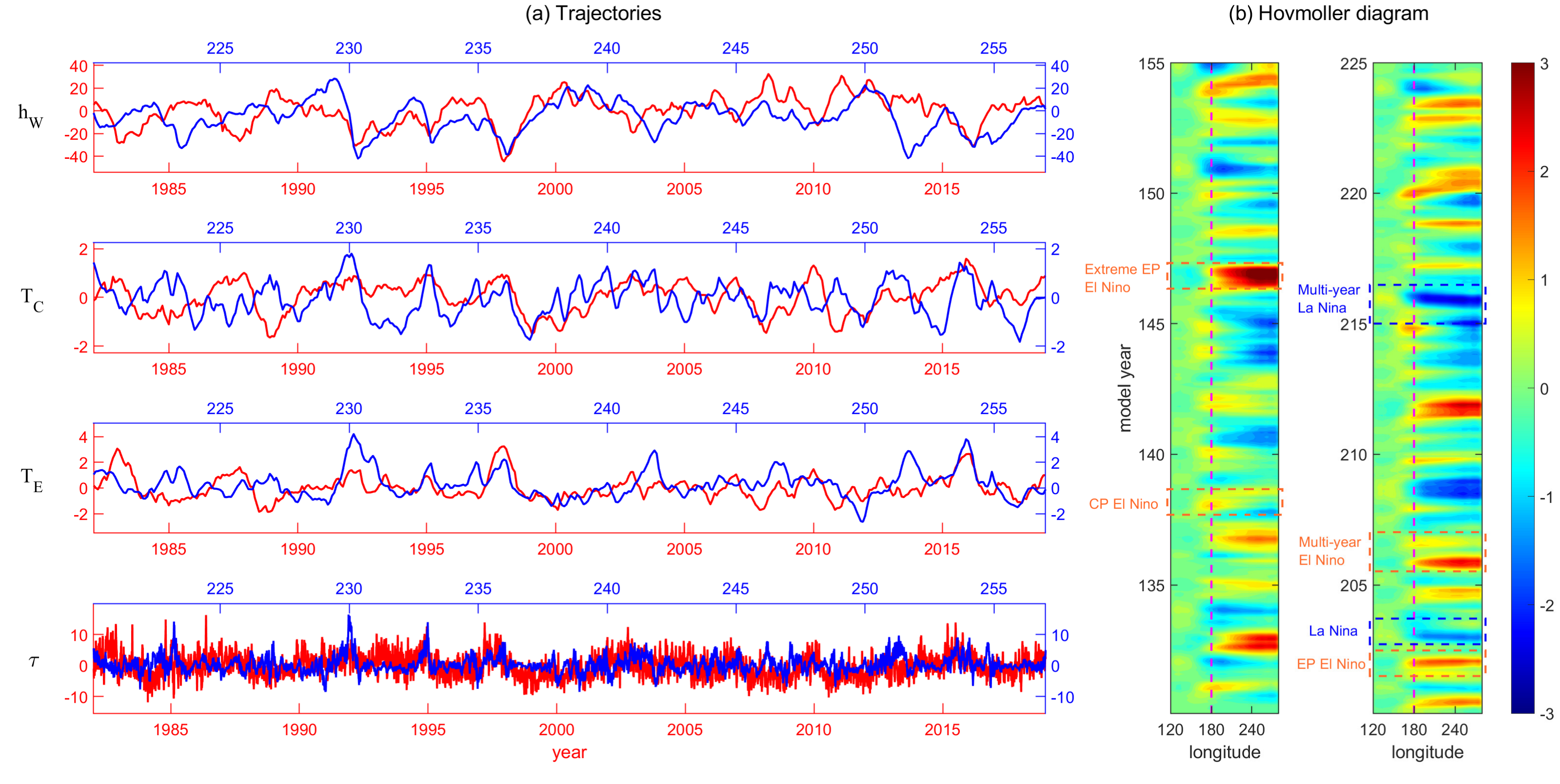}
    \caption{Simulation of the 4-dimensional $I_sI_aM$. Panel (a): time series of the model (blue), which is compared with observations (red) in a qualitative way. Panel (b): Hovmoller diagram of the spatiotemporal pattern of the model simulation. Several different types of ENSO events are marked by dashed boxes, where the orange ones stand for El Ni\~no events and the blue ones represent the La Ni\~no events. }
    \label{Fig:PartII_4d_tra}
\end{figure}

Figure \ref{Fig:PartII_4d_tra} shows the time series and the spatiotemporal pattern of the 4-dimensional $I_sI_aM$. The time series in Panel (a) illustrates a clear relationship between the wind $\tau$ and the EP events in the model, for example, around $t=230$, $236$, and $254$. The coupled relationship between the three interannual variabilities in the model is also quite similar to those in the observations. The spatiotemporal patterns of different types of ENSO events are presented in Panel (b), including EP El Ni\~no ($t=202$), CP El Ni\~no ($t=138$), La Ni\~na ($t=203$), extreme El Ni\~no ($t=147$), multi-year El Ni\~no  ($t=206$ to $t=207$) and multi-year La Ni\~na ($t=215$ to $t=216$).

\subsection{Improving the model performance by incorporating latent variables}\label{Subsec:Latent}

It has been shown that the 3-dimensional interannual model fails to reproduce several crucial features of ENSO diversity. This is because the model lacks certain critical state variables. In general, perfectly pre-determining the most appropriate set of state variables is not always an easy task. Therefore, it is beneficial if the auto-learning framework can identify the missing components in the target model and include additional suitable processes to compensate for the missing information automatically. As the physical meaning of the additional processes is unknown at this stage, the associated variables are treated as latent variables. Once the time series of these latent variables are recovered, physical justifications can be provided to explain the resulting model.

To improve the performance of the 3-dimensional interannual model, one latent variable is added to the system. The auto-learning framework then involves a joint learning procedure to derive the equations for the three pre-determined state variables $(h_W, T_C, T_E)$ and the latent variable from the synthetic data generated from the reference model. Notably, the time series is only available for the three pre-determined state variables, as the meaning of the latent variable remains unclear at this stage. Therefore, the learning procedure here fundamentally differs from deriving the 4-dimensional $I_aI_sM$, in which all the four state variables have observations. An iterative method is applied to alternate between inferring the time series of the latent variables and deriving the model structure \citep{chen2022causality}. The output of this method is a new model involving existing state variables and the latent variable after iteration, as shown in Figure \ref{Fig:Overview_CELA}. For simplicity, only additive noise is assumed in the latent process.

\begin{figure}[ht]
    \hspace*{-0cm}\includegraphics[width=1.0\textwidth]{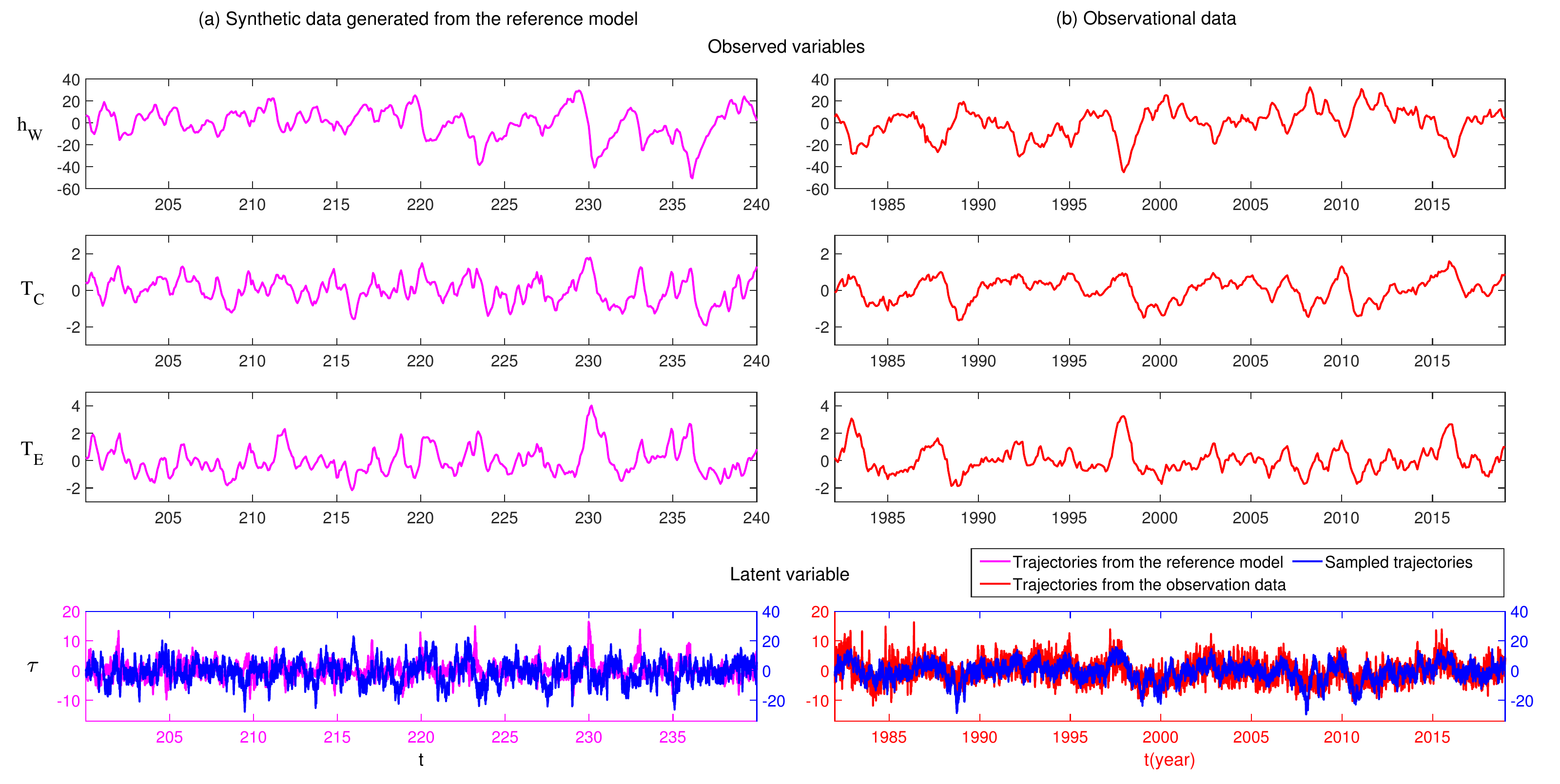}
    \caption{Comparison of observed time series for interannual variables $(h_W, T_C, T_E)$ and a sampled latent variable trajectory derived from the auto-learning framework. Panel (a): the trajectories (pink) of four variables $(h_W, T_C, T_E, \tau)$ from the 6-dimensional reference model and a sampled trajectory (blue) of the latent variable $\tau$ based on the new 4-dimensional model from the learning framework when the time series of $(h_W, T_C, T_E)$ from the reference model is given. Panel (b): Similar to Panel (a) except the time series from the reference model is replaced by the observational data (red).}
    \label{Fig:PartII_ENSO_Latent3_sample}
\end{figure}

To figure out if the latent variable mimics the behavior of a certain unobserved state variable, a comparison between the latent variable and this unobserved state variable in a path-wise way is necessary. Here, synthetic data generated from the reference model and real observations are taken as truth in the following two experiments, respectively. Based on the newly derived 4-dimensional model from the auto-learning framework and true trajectories of existing state variables, the time series of the latent variable could be sampled via the conditional sampling method. Figure \ref{Fig:PartII_ENSO_Latent3_sample} shows the observed time series of the three interannual variables $(h_W, T_C, T_E)$ and one sampled time series of the latent variable using the ensemble Kalman smoother data assimilation approach (see also Section \ref{Subsec:criterion_1a}). Panel (a) shows a sampled trajectory based on the time series of $(h_W, T_C, T_E)$ generated from the reference model, while Panel (b) utilizes the observational data set as the input. A significant finding is that the sampled trajectory of the latent variable has a very similar pattern with the atmospheric wind in both cases. Therefore, the auto-learning framework succeeds in identifying the dominant component missing in the 3-dimensional interannual model. The latent variable thus provides a natural surrogate of the atmospheric wind in this new model. For this reason, the latent variable is denoted by $\tau$. Note that the amplitude of the sampled trajectories differs from that of the actual observed wind burst $\tau$ (see the last row of Figure \ref{Fig:PartII_ENSO_Latent3_sample}). This is not surprising as $\tau$ multiplied by the wind stress coefficient is the actual forcing to the SST processes; for example, the forcing is given by $L_{34} {\tau}$ in \eqref{g_Te}. Neither $L_{34}$ nor ${\tau}$ is uniquely determined but the actual wind forcing $L_{34} {\tau}$ can be learned accurately. If the wind stress coefficient is known, then the latent variable is naturally the wind burst.

\begin{figure}[ht]
    \hspace*{-0cm}\includegraphics[width=1.0\textwidth]{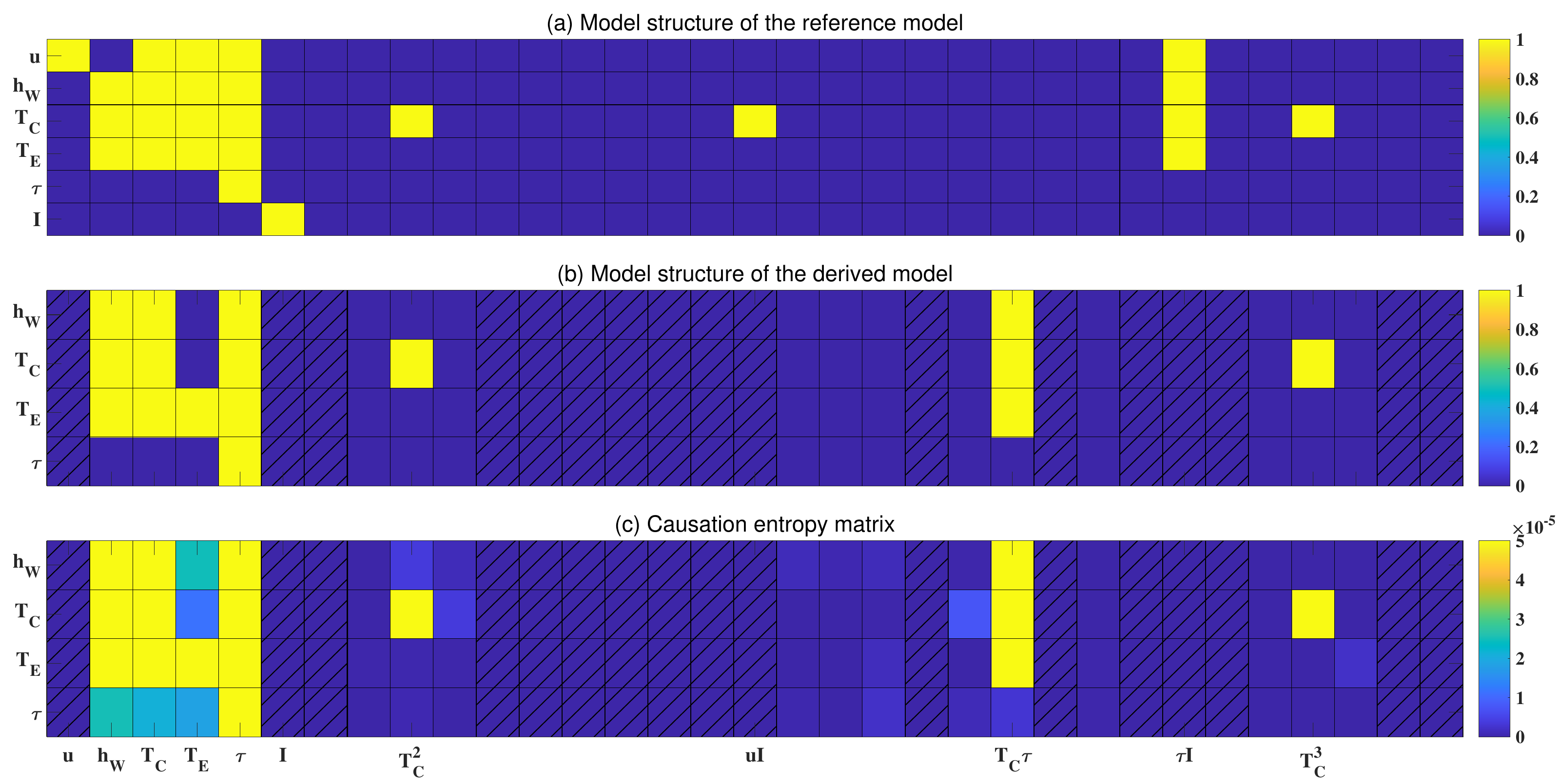}
    \caption{Comparison of the model structures of the 6-dimensional reference model (Panel (a)) and the resulting new 4-dimensional model containing one latent variable $\tau$ (Panel (b)). Panel (c) shows the causation entropy of the new model. The dark blue shadings show the terms not included in the pre-determined function library. }
    \label{Fig:PartII_ENSO_Latent3_C}
\end{figure}

Figure \ref{Fig:PartII_ENSO_Latent3_C} displays the structure of the resulting new 4-dimensional model (Panel (b)) and compares it with the 6-dimensional reference model (Panel (a)). The latent variable is denoted by $\tau$ in this figure since it has been identified as a surrogate for the atmospheric wind. The resulting new 4-dimensional model provides a very similar structure to the reference model in terms of the dynamics $(h_W, T_C, T_E)$, which are dominated by linear terms but contain a cubic damping for the process of $T_C$. It is worth highlighting that $T_C\tau$ appears on the right-hand side of the processes for the three interannual variables. Since only an additive noise is used in the latent process, this $T_C\tau$ term provides an equivalent multiplicative noise, indicating that the effective strength of the wind depends on the SST.

\begin{figure}[ht]
    \hspace*{-0cm}\includegraphics[width=1.0\textwidth]{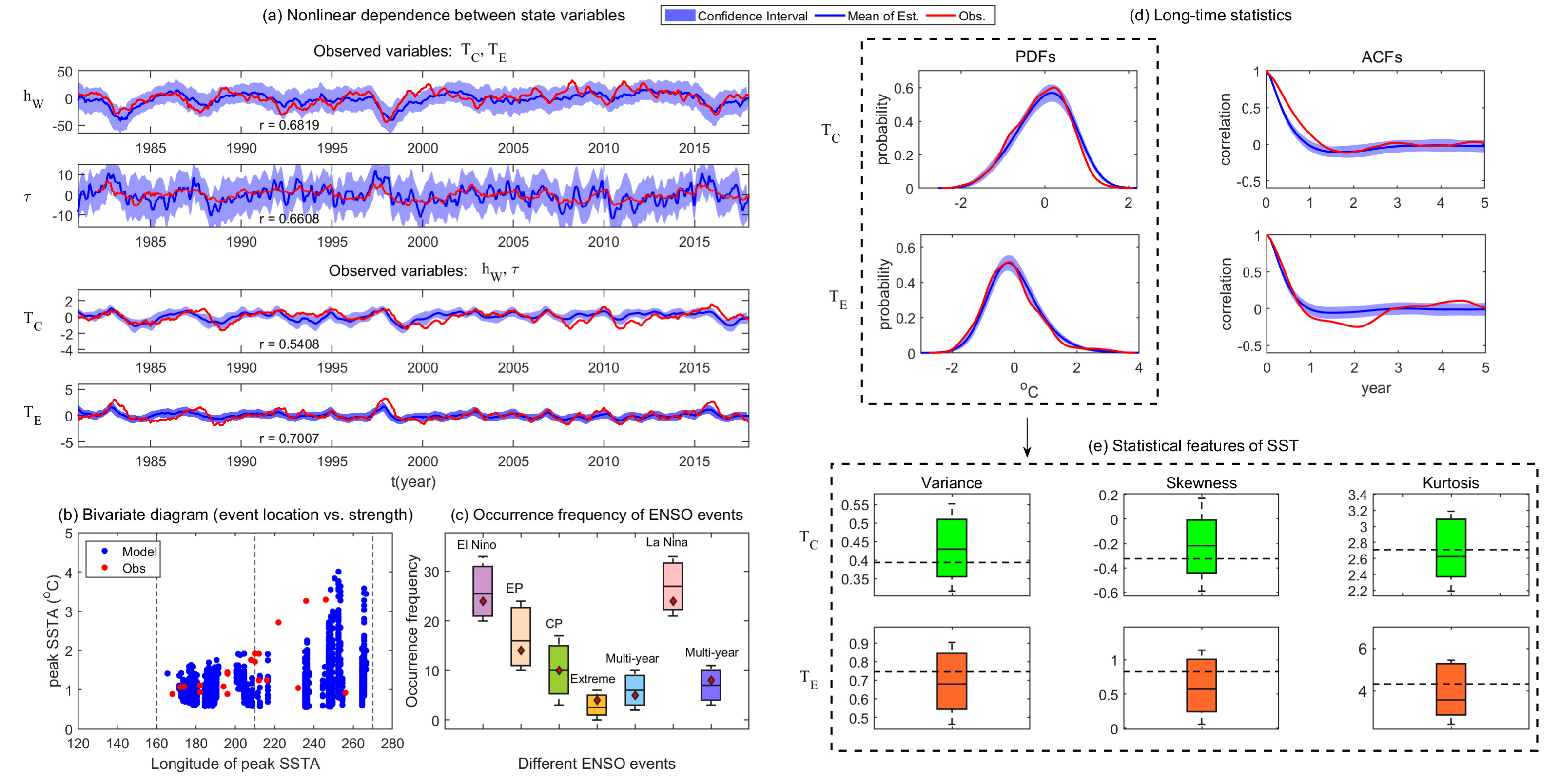}
    \caption{Validation tests for the new 4-dimensional model derived from the auto-learning framework. Panel (a): assessing the fundamental dynamical properties. Note that the wind time series has been plotted as a monthly average for illustration purposes. Panel (b): the bivariate diagram of the locations and strengths of different ENSO events. Panel (c): the occurrence frequency of different ENSO events. In the box plot, the central mark on each box denotes the median, and the bottom and top edges present the 5th and 95th percentiles, respectively. Panel (d): long-time statistics including PDFs and ACFs of $T_C$ and $T_E$. Panel (e): additional supporting information of PDFs in Panel (d). It presents the statistical features of SST, containing the variance, skewness, and kurtosis of $T_C$ and $T_E$.}
    \label{Fig:PartII_ENSO_Latent3_validations}
\end{figure}

Figure \ref{Fig:PartII_ENSO_Latent3_validations} includes the dynamical and statistical features of this new 4-dimensional model. These features are mostly missed by the 3-dimensional interannual model (see Figures \ref{fig:PartII_DA_new}--\ref{Fig:EC_II}) but are overall accurately recovered by the new 4-dimensional model. Panel (a) shows the results of two experiments utilizing Bayesian Inference. In the first experiment, the observed time series of $T_C$, $T_E$ are used to estimate the state of $h_W$ and the latent variable (the surrogate for the wind burst $\tau$ in observational data), where the amplitude of the recovered wind time series has been adjusted to match the observations. In the second one, the observed time series of $h_W$, $\tau$ are used to recover the state of $T_C$ and $T_E$. Among the two experiments, the recovered time series of $h_W$ is more accurate than that of the 3-dimensional model, which indicates that the new 4-dimensional model captures a better oscillation mechanism between SST and thermocline with the contribution of the latent variable. Next, Panel (b) presents the bivariate plot of the location and strength of El Ni\~no events. The new 4-dimensional model captures the feature that the average strength of CP events is lower than those of EP events. Besides, the model can generate a sufficient number of extreme events. This is also supported by the box plot of the occurrence frequency of different ENSO events in Panel (c). The occurrence frequencies of different ENSO events from observational data (the red diamonds) lie in the corresponding box. Similarly, the model performs well in reproducing long-term statistics, especially the PDFs, which can be seen in Panel (d). Recall that a noticeable bias was found in the PDF of $T_C$ from the 3-dimensional model in Panel (c) of Figure \ref{Fig:PartII_Category2}. With the addition of the latent variable, the new 4-dimensional model can accurately capture the non-Gaussian features. Further validation is shown in Panel (e), which contains box plots of the variance, skewness, and kurtosis. It illustrates similar profiles among these aspects with observations (black dashed line). In conclusion, the latent variable serving as the wind process provides essential feedback to the interannual variables. It allows the new model to reproduce the observed features of ENSO diversity.

\section{Further Discussions and Conclusion}\label{Sec:Conclusion}

\subsection{Summary}

In this paper, a physics-informed auto-learning framework is built to derive a hierarchy of stochastic conceptual models with different state variables and degrees of freedom for ENSO diversity. These models are intercompared based on validation criteria. Except for the 3-dimensional interannual model, all other models have comparably skillful results in reproducing many observed features of ENSO diversity. Therefore, the minimum model for ENSO diversity should contain at least three interannual variables describing the thermocline depth in the western Pacific, SST in the EP and the CP, and one intraseasonal variable for wind bursts. Without the latter, the resulting 3-dimensional interannual model suffers from underestimating extreme events and biased PDFs. In addition to the dimensions of the system, the minimum complexity of the model is also discussed. It is also shown that certain nonlinearities are crucial in reproducing the observed non-Gaussian features and realistic amplitudes of different ENSO events. Finally, the auto-learning framework can also recognize the sources of model error and systematically incorporate physically explainable latent variables into the model derivation. These latent variables significantly improve the performance of the resulting models. Therefore, the auto-learning framework offers a powerful tool to derive models when only limited prior knowledge is available.

\subsection{The multiscale nature in modeling ENSO diversity}

It has already been shown in the previous section that the 3-dimensional interannual model fails to characterize several crucial features of ENSO diversity. In contrast, the new 4-dimensional model that includes the intraseasonal atmospheric wind variations significantly improves the resulting dynamics and statistics. This confirms the necessity of having intraseasonal variability in modeling ENSO diversity.

On the other hand, the role of decadal variability in contributing to the ENSO diversity can be understood by comparing the 4-dimensional $I_a I_s M$ with the 5-dimensional $I_a I_s D M$ and the 6-dimensional $I_a I_s D M A$. It turns out that the results based on the validation metrics do not have a significant change when the decadal variable is excluded from the models. However, caution is needed to explain such a result. All the results are based on 70-year segments of data as same as the length of observations, which could be too short for the random modulation of the decadal variability to play a significant role in changing the model statistics. In other words, its role in modifying the statistics may not be distinguishable from certain random noise. It is worth highlighting that decadal variability is essential \citep{dieppois2021enso} in situations with regime switching, such as the climate change scenario \citep{yu2007decadal,yu2013identifying}. This was the guideline for designing the multiscale model in \citep{chen2022multiscale}, where the decadal variability determines the CP- or EP-favored regimes. If the mean state of the decadal variability is changed, then it will affect the occurrence frequency of the CP and EP events \citep{chen2022multiscale}, which cannot be easily reproduced by simple random noises. The model performance for such a sensitivity experiment is not included in the model validation criteria adopted here. These criteria are primarily designed to examine the skill of the models in reproducing the current climate.

\subsection{Applicability to CMIP and paleo databases}

The auto-learning framework utilizes information theory to identify model structures requiring long-time datasets. While real-world observations are limited, model simulations and paleo reconstructions provide ample data for training conceptual models. Coupled Model Intercomparison Project (CMIP) model outputs contain multi-century projections enabling robust structure extraction \citep{eyring2016overview}. Paleo proxy records also have long climate time series that could inform minimal model development for different climate epochs \citep{brown2008mid}. Therefore, applying the auto-learning approach with CMIP and paleo databases is an intriguing future direction. The framework may extract underlying dynamics from complex model projections. Paleo-informed models could also give insights into ENSO evolution under past climate regimes \citep{ford2012deep}. The feasibility of extracting structures from such alternative datasets merits in-depth exploration.

\subsection{Relation to recharge oscillator models}

The recharge oscillator paradigm established in \citep{jin1997equatorial} represented a foundation in ENSO conceptual modeling. Recent studies have sought to extend the recharge oscillator framework by incorporating additional processes and dimensions. The 6-dimensional conceptual model used here to generate synthetic data follows this extended recharge oscillator direction. This work demonstrates that the minimum requirements for a skillful ENSO diversity model should include three interannual states and intraseasonal wind bursts. The results from this work provide helpful guidance on potential pathways to improve recharge oscillator models through appropriate dimensionality expansion. One guideline for advancing recharge oscillator models includes incorporating diversity in ENSO spatial patterns \citep{geng2022enso, chen2023simple}. Thus, the conceptual modeling approach furnishes valuable tools and insights to guide recharge oscillator model development toward more comprehensive ENSO simulation.

\subsection{Potential future directions}

In addition to guiding the development of more sophisticated models directly capturing the spatial patterns of ENSO diversity, the set of stochastic conceptual models developed in this work has several potential applications. These models can be used for multi-model data assimilation and forecast for the large-scale features of ENSO diversity. In particular, the stochastic nature of these models allows statistical tools to analyze the predictability and sensitivity analysis \citep{fang2023quantifying, chen2023stochastic}.

\section*{Acknowledgments}
The research of N.C. is partially funded by ARO W911NF-23-1-0118. Y.Z. is partially supported by this grant.
The research of X.F. is supported by Guangdong Major Project of Basic and Applied Basic Research (Grant No. 2020B0301030004), the Ministry of Science and Technology of the People's Republic of China (Grant No. 2020YFA0608802) and the National Natural Science Foundation of China (Grant Nos. 42192564 and 41805045). The research of J.V. was supported by the French Agence Nationale pour la Recherche (ANR) ARiSE grant (ANR-18-CE01-0012). The authors thank Dr. Fei-Fei Jin for useful discussions. 

\clearpage
\section{Appendix: The Details of the Stochastic Conceptual Models and Parameters Resulting from the Auto-Learning Framework}
\subsection{The interannual intraseasonal decadal model with advection}
The interannual intraseasonal decadal model with advection ($I_a I_s D M A$) is given as follows:
\begin{subequations}
    \begin{align}
    \frac{\d u}{\d t}= & L_{11}  u+L_{12}  \tau+L_{13} +\sigma_u \dot{W}_u, \label{a_u}\\
    \frac{\d h_{W}}{\d t}= & L_{21}  h_{W}+L_{22}  T_{E}+L_{23}  \tau+L_{24}  \tau I+L_{25} +\sigma_h  \dot{W}_h, \label{a_h}\\
    \frac{\d T_{C}}{\d t}= & L_{31}  h_{W}+L_{32}(t) T_{C}+L_{33}  T_{E}+L_{34} \tau+L_{35}(t)  T_{C}^2+L_{36}  u I+L_{37}  \tau I \notag\\
    & +L_{38}(t)  T_{C}^3+L_{39} + \sigma_{C}  \dot{W}_C, \label{a_Tc}\\
    \frac{\d T_{E}}{\d t}= & L_{41}  h_{W}+L_{42} T_{C}+L_{43}(t) T_{E}+L_{44}  \tau+L_{45} \tau I+L_{46} +\sigma_E  \dot{W}_E, \label{a_Te}\\
    \frac{\d \tau}{\d t}= & L_{51}  \tau+L_{52} +\sigma_\tau \left(T_{C}\right) \dot{W}_{\mathrm{\tau}}, \label{a_tau}\\
    \frac{\d I}{\d t}= & L_{61}  I +L_{62} +\sigma_I (I) \dot{W}_I, \label{a_I}
\end{align}
\end{subequations}
where the parameters are listed in Table \ref{Table:Model_a}.
\begin{table}[ht]
\centering
\begin{tabular}{llllll}
 \hline
$L_{11}$   & $-0.1400$  & $L_{31}$     & $0.4540$                                & $L_{41}$ & $0.4493$                                                                                           \\
$L_{12}$   & $-0.0428$ & $L_{32}$     & $-0.4505-0.2989s_1$                   & $L_{42}$ & $-0.2830$                                                                                          \\
$L_{13}$   & $0.0000$  & $L_{33}$     & $0.2850$                                & $L_{43}$ & $-0.0558-0.3618s_2-0.1788s_3$ \\
$L_{21}$   & $-0.1663$ & $L_{34}$     & $0.1983$                                & $L_{44}$ & $0.2470$                                                                                           \\
$L_{22}$   & $-0.0694$ & $L_{35}$     & $-3.1142-1.2090s_1$                  & $L_{45}$ & $-0.0245$                                                                                          \\
$L_{23}$   & $-0.1007$ & $L_{36}$     & $0.1218$                                & $L_{46}$ & $0.0001$                                                                                           \\
$L_{24}$   & $0.0100$  & $L_{37}$     & $-0.0196$                               & $L_{51}$ & $-1.9942$                                                                                          \\
$L_{25}$   & $0.0001$  & $L_{38}$     & $-15.6559-6.2024s_1$                  & $L_{52}$ & $0.0045$                                                                                           \\
$\sigma_u$ & $0.0310$  & $L_{39}$     & $0.0177$                                & $L_{61}$ & $-0.0323$                                                                                          \\
$\sigma_h$ & $0.0155$  & $s_1$        & $\sin \left(2\pi t/6\right)$             & $L_{62}$ & $0.0639$                                                                                           \\
$\sigma_C$ & $0.0310$  & $s_2$        & $\sin(2\pi t/6 +2\pi/6)$ & $s_3$    & $\sin(2\pi t/3+2\pi/6)$                                                           \\
$\sigma_E$ & $0.0232$  & $\sigma_{\tau}$ & \multicolumn{3}{l}{$0.8999\left[\tanh \left(4.5 T_C\right)+1\right]\left[1+0.25 \cos \left(2 \pi/6 t\right)\right]$}                       \\ \hline
\end{tabular}\caption{Parameters of the interannual intraseasonal decadal model with advection ($I_a I_s D M A$).} Note that the functions of multiplicative noises are not included in the library since they could not be identified from the framework. The multiplicative noise of $\tau$ and $I$ are assumed to be known and they are the same as the reference model.\label{Table:Model_a}
\end{table}

\subsection{The interannual intraseasonal model with advection}

The interannual intraseasonal model with advection ($I_a I_s M A$) is given as follows:
\begin{subequations}
\begin{align}
\frac{\d u}{\d t}= & L_{11}  u+L_{12}   \tau+L_{13} +\sigma_u  \dot{W}_u, \label{b_u}\\
\frac{\d h_{W}}{\d t}= & L_{21}  h_{W}+L_{22} T_C + L_{23} T_{E}+L_{24}  \tau+L_{25} +\sigma_h  \dot{W}_h, \label{b_h}\\
\frac{\d T_{C}}{\d t}= & L_{31}  u+L_{32}  h_{W}+L_{33}(t)  T_{C}+L_{34}  T_{E}+L_{35}  \tau + L_{36}(t)  T_{C}^2 + L_{37}(t)  T_{C}^3 \notag\\
& +L_{38}+\sigma_C  \dot{W}_C, \label{b_c} \\
\frac{\d T_{E}}{\d t}= & L_{41}  h_{W}+L_{42}  T_{C}+L_{43}(t)  T_{E}+L_{44}  \tau+L_{45} +\sigma_E  \dot{W}_E, \label{b_e}\\
\frac{\d \tau}{\d t}= & L_{51} \tau+L_{52} +\sigma_\tau \left(T_{C}\right) \dot{W}_\tau, \label{b_tau}
\end{align}
\end{subequations}
where the parameters are listed in Table \ref{Table:Model_b}.
\begin{table}[ht]
\centering
\begin{tabular}{llllll}
 \hline
$L_{11}$   & $-0.1400$ & $L_{31}$     & $0.2183$                                               & $L_{41}$               & $0.4443$                                        \\
$L_{12}$   & $-0.0428$ & $L_{32}$     & $0.4222$                                               & $L_{42}$               & $-0.2741$                                       \\
$L_{13}$   & $0.0000$  & $L_{33}$     & $-0.4061-0.2989s_1$                                  & $L_{43}$               & $-0.0553-0.3597s_2-0.1791s_3$               \\
$L_{21}$   & $-0.1481$ & $L_{34}$     & $0.2455$                                               & $L_{44}$               & $0.1976$                                        \\
$L_{22}$   & $-0.0501$ & $L_{35}$     & $0.1568$                                               & $L_{45}$               & $0.0001$                                        \\
$L_{23}$   & $-0.0456$ & $L_{36}$     & $-2.8721-1.1140s_1$                                 & $L_{51}$               & $-1.9942$                                       \\
$L_{24}$   & $-0.0793$ & $L_{37}$     & $-14.5992-5.6033s_1$                                 & $L_{52}$               & $0.0045$                                       \\
$L_{25}$   & $0.0000$  & $L_{38}$     & $0.0166$                                               &                &                                       \\
$\sigma_u$ & $0.0310$  & $s_1$        & $\sin \left(2\pi t/6\right)$                  &                        &                                               \\
$\sigma_h$ & $0.0155$  & $s_2$        & $\sin(2\pi t/6+2\pi/6)$                &                        &                                               \\
$\sigma_C$ & $0.0310$  & $s_3$        & $\sin(2\pi t/3+2\pi/6)$               &                        &                                               \\
$\sigma_E$ & $0.0232$  & $\sigma_{\tau}$ & \multicolumn{3}{l}{$0.8999\left[\tanh \left(4.5 T_C\right)+1\right]\left[1+0.25 \cos \left(2 \pi t/6\right)\right]$} \\ \hline
\end{tabular}\caption{Parameters of the interannual intraseasonal model with advection ($I_a I_s M A$).}\label{Table:Model_b}
\end{table}

\subsection{The interannual intraseasonal decadal model}
The interannual intraseasonal decadal model ($I_a I_s D M$) is given as follows:
\begin{subequations}
\begin{align}
\frac{\d h_{W}}{ \d t}= & L_{11}  h_{W}+L_{12} T_E +  L_{13} \tau +L_{14}  \tau  I+L_{15} +\sigma_h  \dot{W}_h, \label{c_h}\\
\frac{\d T_{C}}{ \d t}= & L_{21}  h_{W}+L_{22}(t)  T_{C}+L_{23}T_{E}+L_{24}   \tau+L_{25}(t)  T_{C}^2 + L_{26} T_C I + L_{27} T_E I + L_{28} \tau  I \notag\\
& +L_{29}(t)  T_{C}^3  +L_{2,10} +\sigma_C  \dot{W}_C, \label{c_Tc}\\
\frac{\d T_{E}}{ \d t}= & L_{31}  h_{W}+L_{32}  T_{C}+L_{33}(t)  T_{E}+L_{34}   \tau+L_{35}  \tau  I+L_{36} + \sigma_E  \dot{W}_E, \label{c_Te}\\
\frac{\d \tau}{ \d t}= & L_{41}   \tau+L_{42} +\sigma_\tau \left(T_{C}\right) \dot{W}_{\tau}, \label{c_tau}\\
\frac{\d I}{\d t}= & L_{51}   I+L_{52} +\sigma_I (I) \dot{W}_I, \label{c_I}
\end{align}
\end{subequations}
where the parameters are listed in Table \ref{Table:Model_c}.

\begin{table}[ht]
\centering
\begin{tabular}{llllll}
\hline
$L_{11}$ & $-0.1663$             & $L_{27}$       & $-0.0563$                                          & $L_{41}$                     & $-1.9942$                                     \\
$L_{12}$ & $-0.0694$             & $L_{28}$       & $-0.0203$                                          & $L_{42}$                     & $0.0045$                                      \\
$L_{13}$ & $-0.1007$             & $L_{29}$       & $-15.3784-6.2253s_1$                             & $L_{51}$                     & $-0.0323$                                     \\
$L_{14}$ & $0.0100$              & $L_{2,10}$     & $0.0166$                                           & $L_{52}$                     & $0.0639$                                      \\
$L_{15}$ & $0.0001$              & $L_{31}$       & $0.4493$                                           & $s_1$                        & $\sin \left(2\pi t/6\right)$                 \\
$L_{21}$ & $0.5030$              & $L_{32}$       & $-0.2830$                                          & $s_2$                        & $\sin(2\pi t/6+2\pi/6)$                       \\
$L_{22}$ & $-0.4785-0.2993s_1$ & $L_{33}$       & $-0.0558-0.3518s_2-0.1788s_3$                  & $s_3$                        & $\sin(2\pi t/3+2\pi/6)$                      \\
$L_{23}$ & $0.3142$              & $L_{34}$       & $0.2470$                                           & $\sigma_{h}$                 & $0.0155$                                      \\
$L_{24}$ & $0.1981$              & $L_{35}$       & $-0.0245$                                          & $\sigma_{C}$                 & $0.0310$                                      \\
$L_{25}$ & $-2.9602-1.1545s_1$ & $L_{36}$       & $0.0001$                                           & $\sigma_E$                   & $0.0232$                                      \\
$L_{26}$ & $0.0555$              & $\sigma_{\tau}$ & \multicolumn{3}{l}{$0.8999 \left[\tanh \left(4.5 T_C\right)+1\right]\left[1+0.25 \cos \left(2 \pi t/6\right)\right]$} \\ \hline
\end{tabular}\caption{Parameters of the interannual intraseasonal decadal model ($I_a I_s D M$).}\label{Table:Model_c}
\end{table}

\subsection{The interannual intraseasonal model}
The interannual intraseasonal model ($I_a I_s M$) is given by
\begin{subequations}\label{IaIsM}
\begin{align}
\frac{ \d h_{W}}{ \d t}= & L_{11}  h_{W}+L_{12} T_C + L_{13} T_{E}+L_{14} \tau+L_{15} +\sigma_h  \dot{W}_h, \label{d_h}\\
\frac{ \d T_{C}}{ \d t}= & L_{21}  h_{W}+L_{22}(t)  T_{C}+L_{23}  T_{E}+L_{24}   \tau+L_{25}(t)  T_{C}^2+L_{26}(t)  T_{C}^3 \notag\\
& +L_{27} +\sigma_C  \dot{W}_C, \label{d_Tc} \\
\frac{ \d T_{E}}{ \d t}= & L_{31}  h_{W}+L_{32}  T_{C}+L_{33}(t)  T_{E}+L_{34}   \tau+L_{35} +\sigma_E  \dot{W}_E, \label{d_Te}\\
\frac{ \d \tau}{ \d t}= & L_{41}   \tau+L_{42} +\sigma_\tau \left(T_{C}\right) \dot{W}_{\tau}, \label{d_tau}
\end{align}
\end{subequations}
where the parameters are listed in Table \ref{Table:Model_d}.

\begin{table}[ht]
\centering
\begin{tabular}{llllll}
\hline
$L_{11}$       & $-0.1481$                                  & $L_{25}$                     & $-2.7519-1.0865s_1$                                 & $L_{41}$     & $-1.9942$                     \\
$L_{12}$       & $-0.0501$                                  & $L_{26}$                     & $-14.2956-5.7679s_1$                                & $L_{42}$     & $0.0045$                      \\
$L_{13}$       & $-0.0456$                                  & $L_{27}$                     & $0.0158$                                              & $s_1$        & $\sin \left(2\pi t/6\right)$ \\
$L_{14}$       & $-0.0793$                                  & $L_{31}$                     & $0.4446$                                              & $s_2$        & $\sin(2\pi t/6+2\pi/6)$       \\
$L_{15}$       & $0.0000$                                   & $L_{32}$                     & $-0.2741$                                             & $s_3$        & $\sin(2\pi t/3+2\pi/6)$      \\
$L_{21}$       & $0.4861$                                   & $L_{33}$                     & $-0.0553-0.3597s_2-0.1791s_3$                     & $\sigma_{h}$ & $0.0155$                      \\
$L_{22}$       & $-0.3443-0.3014s_1$                      & $L_{34}$                     & $0.1976$                                              & $\sigma_{C}$ & $0.0310$                      \\
$L_{23}$       & $0.1813$                                   & $L_{35}$                     & $0.0000$                                              & $\sigma_E$   & $0.0232$                      \\
$L_{24}$       & $0.1558$                                   &   $\sigma_{\tau}$ & \multicolumn{3}{l}{$0.8999 \left[\tanh \left(4.5 T_C\right)+1\right]\left[1+0.25 \cos \left(2 \pi t/6\right)\right]$}                            \\ \hline
\end{tabular}\caption{Parameters of the interannual intraseasonal model ($I_a I_s M$).}\label{Table:Model_d}
\end{table}

\subsection{The interannual model}
The interannual model ($I_a M$) is given by
\begin{subequations}
    \begin{align}
    \frac{\d h_{W}}{\d t}= & L_{11} h_{W}+L_{12} T_{C}+L_{13} T_{E}+L_{14} T_{C}^2+L_{15}+\sigma_h \dot{W}_h, \label{e_h}\\
    \frac{\d T_{C}}{\d t}= & L_{21} h_{W}+L_{22}(t) T_{C}+L_{23} T_{E}+L_{24} \mathrm{h}_{W} T_{C}+L_{25} T_{C} T_{E}+L_{26} +\sigma_C \dot{W}_C, \label{e_Tc}\\
    \frac{\d T_{E}}{\d t}= & L_{31} h_{W}+ L_{32}(t) T_{E}+L_{33} T_{C}^2 + L_{34}(t) T_E^2 + L_{35}(t) T_{E} ^3 + L_{36}+\sigma_E \dot{W}_E, \label{e_Te}
    \end{align}
\end{subequations}
where the  parameters are listed in Table \ref{Table:Model_e}.

\begin{table}[ht]
\centering
\begin{tabular}{llllll}
\hline
$L_{11}$ & $-0.0678$             & $L_{24}$ & $0.6729$                          & $L_{36}$     & $-0.0078$                     \\
$L_{12}$ & $-0.1927$             & $L_{25}$ & $-1.1864$                         & $s_1$        & $\sin \left(2\pi t/6\right)$ \\
$L_{13}$ & $-0.0593$             & $L_{26}$ & $0.0070$                          & $s_2$        & $\sin(2\pi t/6+2\pi/6)$       \\
$L_{14}$ & $-0.6729$             & $L_{31}$ & $0.2719$                          & $s_3$        & $\sin(2\pi t/3+2\pi/6)$      \\
$L_{15}$ & $0.0049$              & $L_{32}$ & $0.0523-0.2549s_2-0.1792s_3$  & $\sigma_{h}$ & $0.0156$                      \\
$L_{21}$ & $0.3146$              & $L_{33}$ & $1.1864$                          & $\sigma_{C}$ & $0.0311$                      \\
$L_{22}$ & $-0.2553-0.1799s_1$ & $L_{34}$ & $-0.0134+0.0284s_2+0.0381s_3$ & $\sigma_{E}$ & $0.0236$                      \\
$L_{23}$ & $0.2340$              & $L_{35}$ & $-0.2611-0.3158s_2+0.1207s_3$ &              &                             \\ \hline
\end{tabular}\caption{Parameters of the interannual model ($I_a M$).}\label{Table:Model_e}
\end{table}

\subsection{The linear model}
The linear 6-dimensional model is given as follows:
\begin{subequations}\label{model_linear}
    \begin{align}
    \frac{\d u}{\d t}= & L_{11}  u+L_{12}  \tau+L_{13} +\sigma_u \dot{W}_u, \label{f_u}\\
    \frac{\d h_{W}}{\d t}= & L_{21}  h_{W}+L_{22}  T_{E}+L_{23}  \tau+L_{24}  \tau I+L_{25} +\sigma_h  \dot{W}_h, \label{f_h}\\
    \frac{\d T_{C}}{\d t}= & L_{31}  h_{W}+L_{32}(t) T_{C}+L_{33}  T_{E}+L_{34} \tau+L_{35}  u I+L_{36}  \tau I \notag\\
    & +L_{37} + \sigma_{C}  \dot{W}_C, \label{f_Tc}\\
    \frac{\d T_{E}}{\d t}= & L_{41}  h_{W}+L_{42} T_{C}+L_{43}(t) T_{E}+L_{44}  \tau+L_{45} \tau I+L_{46} +\sigma_E  \dot{W}_E, \label{f_Te}\\
    \frac{\d \tau}{\d t}= & L_{51}  \tau+L_{52} +\sigma_\tau \left(T_{C}\right) \dot{W}_{\mathrm{\tau}}, \label{f_tau}\\
    \frac{\d I}{\d t}= & L_{61}  I +L_{62} +\sigma_I (I) \dot{W}_I, \label{f_I}
\end{align}
\end{subequations}
where the parameters are listed in Table \ref{Table:Model_linear}.
\begin{table}[ht]
\centering
\begin{tabular}{llllll}
 \hline
$L_{11}$   & $-0.1400$  & $L_{31}$     & $0.4461$                                & $L_{41}$ & $0.4493$                                                                                           \\
$L_{12}$   & $-0.0428$ & $L_{32}$     & $-0.5564-0.2799s_1$                   & $L_{42}$ & $-0.2830$                                                                                          \\
$L_{13}$   & $0.0000$  & $L_{33}$     & $0.1990$                                & $L_{43}$ & $-0.0558-0.3618s_2-0.1788s_3$ \\
$L_{21}$   & $-0.1663$ & $L_{34}$     & $0.1766$                                & $L_{44}$ & $0.2470$                                                                                           \\
$L_{22}$   & $-0.0694$ & $L_{35}$     & $0.1000$                                & $L_{45}$ & $-0.0245$                                                                                          \\
$L_{23}$   & $-0.1007$ & $L_{36}$     & $-0.0162$                                & $L_{46}$ & $0.0001$                                                                                           \\
$L_{24}$   & $0.0100$  & $L_{37}$     & $-0.0007$                               & $L_{51}$ & $-1.9942$                                                                                          \\
$L_{25}$   & $0.0001$  & $s_1$     &   $\sin \left(2\pi t/6\right)$              & $L_{52}$ & $0.0045$                                                                                           \\
$\sigma_u$ & $0.0310$  & $s_2$        & $\sin(2\pi t/6 +2\pi/6)$                    & $L_{61}$ & $-0.0323$                                                                                          \\
$\sigma_h$ & $0.0155$  & $s_3$        & $\sin(2\pi t/3+2\pi/6)$                    & $L_{62}$ & $0.0639$                                                                                           \\
$\sigma_C$ & $0.0310$  & $\sigma_{\tau}$  & \multicolumn{3}{l}{0.8999$ \left[\tanh \left(4.5 T_C\right)+1\right]\left[1+0.25 \cos \left(2 \pi/6 t\right)\right]$}\\
$\sigma_E$ & $0.0232$  &  \\ \hline
\end{tabular}\caption{Parameters of the linear 6-dimensional model.}\label{Table:Model_linear}
\end{table}

\subsection{The new 4-dimensional model with the latent variable}
The new 4-dimensional model with the latent variable is given by
\begin{subequations}\label{new_4d}
\begin{align}
\frac{ \d h_{W}}{ \d t}= & L_{11}  h_{W}+L_{12} T_C + L_{13} \hat{\tau} + L_{14} T_C\hat{\tau} + L_{15} +\sigma_h  \dot{W}_h, \label{g_h}\\
\frac{ \d T_{C}}{ \d t}= & L_{21}  h_{W}+L_{22}(t)T_{C}+L_{23} \hat{\tau}+L_{24}(t)  T_{C}^2 +L_{25} T_C\hat{\tau}+L_{26}(t)  T_{C}^3 \notag\\
& +L_{27} +\sigma_C  \dot{W}_C, \label{g_Tc} \\
\frac{ \d T_{E}}{ \d t}= & L_{31}  h_{W}+L_{32}  T_{C}+L_{33}(t)  T_{E}+L_{34}   \hat{\tau}+L_{35}  T_{C} \hat{\tau}+L_{36} +\sigma_E  \dot{W}_E, \label{g_Te}\\
\frac{ \d \hat{\tau}}{ \d t}= & L_{41} \hat{\tau}+L_{42} +\sigma_{\hat{\tau}}\dot{W}_{\hat{\tau}}, \label{g_tau}
\end{align}
\end{subequations}
where the $\hat{\tau}$ denotes the latent variable, the parameters are listed in Table \ref{Table:new_4d}.

\begin{table}[ht]
\centering
\begin{tabular}{llllll}
\hline
$L_{11}$       & $-0.1003$                                  & $L_{25}$                     & $0.1834$                                 & $L_{41}$     & $-1.5815$                     \\
$L_{12}$       & $-0.0819$                                  & $L_{26}$                     & $-10.1500-1.8734s_1$                                & $L_{42}$     & $-0.1297$                      \\
$L_{13}$       & $-0.0289$                                  & $L_{27}$                     & $0.0122$                                              & $s_1$        & $\sin \left(2\pi t/6\right)$ \\
$L_{14}$       & $-0.1069$                                  & $L_{31}$                     & $0.3222$                                              & $s_2$        & $\sin(2\pi t/6+2\pi/6)$       \\
$L_{15}$       & $0.0012$                                   & $L_{32}$                     & $-0.1863$                                             & $s_3$        & $\sin(2\pi t/3+2\pi/6)$      \\
$L_{21}$       & $0.3470$                                   & $L_{33}$                     & $-0.2110-0.1402s_2-0.0768s_3$                     & $\sigma_{h}$ & $0.0155$                      \\
$L_{22}$       & $-0.2554-0.2584s_1$                        & $L_{34}$                     & $0.0782$                                              & $\sigma_{C}$ & $0.0310$                      \\
$L_{23}$       & $0.0589$                                   & $L_{35}$                     & $0.2580$                                              & $\sigma_E$   & $0.0232$                      \\
$L_{24}$       & $-2.1976-0.5077s_1$                        & $L_{36}$                     & $-0.0017$                                             & $\sigma_{\hat{\tau}}$  &  $2.2034$                           \\
\hline
\end{tabular}\caption{Parameters of the new 4-dimensional model with one latent variable.}\label{Table:new_4d}
\end{table}

\bibliography{references}

\begin{thebibliography}{94}
\providecommand{\natexlab}[1]{#1}
\providecommand{\url}[1]{\texttt{#1}}
\expandafter\ifx\csname urlstyle\endcsname\relax
  \providecommand{\doi}[1]{doi: #1}\else
  \providecommand{\doi}{doi: \begingroup \urlstyle{rm}\Url}\fi

\bibitem[Philander(1983)]{philander1983nino}
S~George~H Philander.
\newblock {E}l {N}i{\~n}o {S}outhern {O}scillation phenomena.
\newblock \emph{Nature}, 302\penalty0 (5906):\penalty0 295--301, 1983.

\bibitem[Ropelewski and Halpert(1987)]{ropelewski1987global}
Chester~F Ropelewski and Michael~S Halpert.
\newblock Global and regional scale precipitation patterns associated with the {E}l {N}i{\~n}o/{S}outhern {O}scillation.
\newblock \emph{Monthly Weather Review}, 115\penalty0 (8):\penalty0 1606--1626, 1987.

\bibitem[Klein et~al.(1999)Klein, Soden, and Lau]{klein1999remote}
Stephen~A Klein, Brian~J Soden, and Ngar-Cheung Lau.
\newblock Remote sea surface temperature variations during {ENSO}: Evidence for a tropical atmospheric bridge.
\newblock \emph{Journal of Climate}, 12\penalty0 (4):\penalty0 917--932, 1999.

\bibitem[McPhaden et~al.(2006)McPhaden, Zebiak, and Glantz]{mcphaden2006enso}
Michael~J McPhaden, Stephen~E Zebiak, and Michael~H Glantz.
\newblock {ENSO} as an integrating concept in earth science.
\newblock \emph{Science}, 314\penalty0 (5806):\penalty0 1740--1745, 2006.

\bibitem[Dai and Wigley(2000)]{dai2000global}
Aiguo Dai and TML Wigley.
\newblock Global patterns of {ENSO}-induced precipitation.
\newblock \emph{Geophysical Research Letters}, 27\penalty0 (9):\penalty0 1283--1286, 2000.

\bibitem[Ashok et~al.(2007)Ashok, Behera, Rao, Weng, and Yamagata]{ashok2007nino}
Karumuri Ashok, Swadhin~K Behera, Suryachandra~A Rao, Hengyi Weng, and Toshio Yamagata.
\newblock El {N}i{\~n}o {M}odoki and its possible teleconnection.
\newblock \emph{Journal of Geophysical Research: Oceans}, 112\penalty0 (C11), 2007.

\bibitem[Kao and Yu(2009)]{kao2009contrasting}
Hsun-Ying Kao and Jin-Yi Yu.
\newblock Contrasting eastern-{P}acific and central-{P}acific types of {ENSO}.
\newblock \emph{Journal of Climate}, 22\penalty0 (3):\penalty0 615--632, 2009.

\bibitem[Kim et~al.(2012)Kim, Kim, and Yeh]{kim2012statistical}
Jin-Soo Kim, Kwang-Yul Kim, and Sang-Wook Yeh.
\newblock Statistical evidence for the natural variation of the central {P}acific {E}l {N}i{\~n}o.
\newblock \emph{Journal of Geophysical Research: Oceans}, 117\penalty0 (C6), 2012.

\bibitem[Capotondi et~al.(2015)Capotondi, Wittenberg, Newman, Di~Lorenzo, Yu, Braconnot, Cole, Dewitte, Giese, Guilyardi, et~al.]{capotondi2015understanding}
Antonietta Capotondi, Andrew~T Wittenberg, Matthew Newman, Emanuele Di~Lorenzo, Jin-Yi Yu, Pascale Braconnot, Julia Cole, Boris Dewitte, Benjamin Giese, Eric Guilyardi, et~al.
\newblock Understanding {ENSO} diversity.
\newblock \emph{Bulletin of the American Meteorological Society}, 96\penalty0 (6):\penalty0 921--938, 2015.

\bibitem[Larkin and Harrison(2005)]{larkin2005global}
Narasimhan~K Larkin and DE~Harrison.
\newblock Global seasonal temperature and precipitation anomalies during {E}l {N}i{\~n}o autumn and winter.
\newblock \emph{Geophysical Research Letters}, 32\penalty0 (16), 2005.

\bibitem[Yu and Kao(2007)]{yu2007decadal}
Jin-Yi Yu and Hsun-Ying Kao.
\newblock Decadal changes of {ENSO} persistence barrier in {SST} and ocean heat content indices: 1958--2001.
\newblock \emph{Journal of Geophysical Research: Atmospheres}, 112\penalty0 (D13), 2007.

\bibitem[Kug et~al.(2009)Kug, Jin, and An]{kug2009two}
Jong-Seong Kug, Fei-Fei Jin, and Soon-Il An.
\newblock Two types of {E}l {N}i{\~n}o events: cold tongue {E}l {N}i{\~n}o and warm pool {E}l {N}i{\~n}o.
\newblock \emph{Journal of Climate}, 22\penalty0 (6):\penalty0 1499--1515, 2009.

\bibitem[Jin(2022)]{jin2022toward}
Fei-Fei Jin.
\newblock Toward understanding {E}l {N}i{\~n}o southern-oscillation's spatiotemporal pattern diversity.
\newblock \emph{Frontiers in Earth Science}, 10:\penalty0 899139, 2022.

\bibitem[Chen and Cane(2008)]{chen2008nino}
Dake Chen and Mark~A Cane.
\newblock El {N}i{\~n}o prediction and predictability.
\newblock \emph{Journal of Computational Physics}, 227\penalty0 (7):\penalty0 3625--3640, 2008.

\bibitem[Jin et~al.(2008)Jin, Kinter, Wang, Park, Kang, Kirtman, Kug, Kumar, Luo, Schemm, et~al.]{jin2008current}
Emilia~K Jin, James~L Kinter, Bin Wang, C-K Park, I-S Kang, BP~Kirtman, J-S Kug, A~Kumar, J-J Luo, J~Schemm, et~al.
\newblock Current status of {ENSO} prediction skill in coupled ocean--atmosphere models.
\newblock \emph{Climate Dynamics}, 31\penalty0 (6):\penalty0 647--664, 2008.

\bibitem[Barnston et~al.(2012)Barnston, Tippett, L'Heureux, Li, and DeWitt]{barnston2012skill}
Anthony~G Barnston, Michael~K Tippett, Michelle~L L'Heureux, Shuhua Li, and David~G DeWitt.
\newblock Skill of real-time seasonal {ENSO} model predictions during 2002--11: {I}s our capability increasing?
\newblock \emph{Bulletin of the American Meteorological Society}, 93\penalty0 (5):\penalty0 631--651, 2012.

\bibitem[Hu et~al.(2012)Hu, Kumar, Jha, Wang, Huang, and Huang]{hu2012analysis}
Zeng-Zhen Hu, Arun Kumar, Bhaskar Jha, Wanqiu Wang, Bohua Huang, and Boyin Huang.
\newblock An analysis of warm pool and cold tongue {E}l {N}i{\~n}os: {A}ir--sea coupling processes, global influences, and recent trends.
\newblock \emph{Climate Dynamics}, 38\penalty0 (9):\penalty0 2017--2035, 2012.

\bibitem[Zheng et~al.(2014)Zheng, Fang, Yu, and Zhu]{zheng2014asymmetry}
Fei Zheng, Xiang-Hui Fang, Jin-Yi Yu, and Jiang Zhu.
\newblock Asymmetry of the {B}jerknes positive feedback between the two types of {E}l {N}i{\~n}o.
\newblock \emph{Geophysical Research Letters}, 41\penalty0 (21):\penalty0 7651--7657, 2014.

\bibitem[Fang et~al.(2015)Fang, Zheng, and Zhu]{fang2015cloud}
Xiang-Hui Fang, Fei Zheng, and Jiang Zhu.
\newblock The cloud-radiative effect when simulating strength asymmetry in two types of {E}l {N}i{\~n}o events using {CMIP5} models.
\newblock \emph{Journal of Geophysical Research: Oceans}, 120\penalty0 (6):\penalty0 4357--4369, 2015.

\bibitem[Sohn et~al.(2016)Sohn, Tam, and Jeong]{sohn2016strength}
Soo-Jin Sohn, Chi-Yung Tam, and Hye-In Jeong.
\newblock How do the strength and type of {ENSO} affect {SST} predictability in coupled models.
\newblock \emph{Scientific Reports}, 6\penalty0 (1):\penalty0 1--8, 2016.

\bibitem[Santoso et~al.(2019)Santoso, Hendon, Watkins, Power, Dommenget, England, Frankcombe, Holbrook, Holmes, Hope, et~al.]{santoso2019dynamics}
Agus Santoso, Harry Hendon, Andrew Watkins, Scott Power, Dietmar Dommenget, Matthew~H England, Leela Frankcombe, Neil~J Holbrook, Ryan Holmes, Pandora Hope, et~al.
\newblock Dynamics and predictability of {E}l {N}i{\~n}o--{S}outhern {O}scillation: an {A}ustralian perspective on progress and challenges.
\newblock \emph{Bulletin of the American Meteorological Society}, 100\penalty0 (3):\penalty0 403--420, 2019.

\bibitem[Stuecker et~al.(2013)Stuecker, Timmermann, Jin, McGregor, and Ren]{stuecker2013combination}
Malte~F Stuecker, Axel Timmermann, Fei-Fei Jin, Shayne McGregor, and Hong-Li Ren.
\newblock A combination mode of the annual cycle and the {E}l {N}i{\~n}o/southern oscillation.
\newblock \emph{Nature Geoscience}, 6\penalty0 (7):\penalty0 540--544, 2013.

\bibitem[Timmermann et~al.(2018)Timmermann, An, Kug, Jin, Cai, Capotondi, Cobb, Lengaigne, McPhaden, Stuecker, et~al.]{timmermann2018nino}
Axel Timmermann, Soon-Il An, Jong-Seong Kug, Fei-Fei Jin, Wenju Cai, Antonietta Capotondi, Kim~M Cobb, Matthieu Lengaigne, Michael~J McPhaden, Malte~F Stuecker, et~al.
\newblock {E}l {N}i{\~n}o-{S}outhern {O}scillation complexity.
\newblock \emph{Nature}, 559\penalty0 (7715):\penalty0 535--545, 2018.

\bibitem[Jin(1997{\natexlab{a}})]{jin1997equatorial}
Fei-Fei Jin.
\newblock An equatorial ocean recharge paradigm for {ENSO}. {P}art {I}: {C}onceptual model.
\newblock \emph{Journal of the Atmospheric Sciences}, 54\penalty0 (7):\penalty0 811--829, 1997{\natexlab{a}}.

\bibitem[Jin(1997{\natexlab{b}})]{jin1997equatorial2}
Fei-Fei Jin.
\newblock An equatorial ocean recharge paradigm for {ENSO}. {P}art {II}: {A} stripped-down coupled model.
\newblock \emph{Journal of the Atmospheric Sciences}, 54\penalty0 (7):\penalty0 830--847, 1997{\natexlab{b}}.

\bibitem[Suarez and Schopf(1988)]{suarez1988delayed}
Max~J Suarez and Paul~S Schopf.
\newblock A delayed action oscillator for {ENSO}.
\newblock \emph{Journal of Atmospheric Sciences}, 45\penalty0 (21):\penalty0 3283--3287, 1988.

\bibitem[Battisti and Hirst(1989)]{battisti1989interannual}
David~S Battisti and Anthony~C Hirst.
\newblock Interannual variability in a tropical atmosphere--ocean model: Influence of the basic state, ocean geometry and nonlinearity.
\newblock \emph{Journal of the Atmospheric Sciences}, 46\penalty0 (12):\penalty0 1687--1712, 1989.

\bibitem[Weisberg and Wang(1997)]{weisberg1997western}
Robert~H Weisberg and Chunzai Wang.
\newblock A western {P}acific oscillator paradigm for the {E}l {N}i{\~n}o-{S}outhern {O}scillation.
\newblock \emph{Geophysical Research Letters}, 24\penalty0 (7):\penalty0 779--782, 1997.

\bibitem[Picaut et~al.(1997)Picaut, Masia, and Du~Penhoat]{picaut1997advective}
Jo{\"e}l Picaut, Fran{\c{c}}ois Masia, and Y~Du~Penhoat.
\newblock An advective-reflective conceptual model for the oscillatory nature of the {ENSO}.
\newblock \emph{Science}, 277\penalty0 (5326):\penalty0 663--666, 1997.

\bibitem[Wang(2001)]{wang2001unified}
Chunzai Wang.
\newblock A unified oscillator model for the el ni{\~n}o--southern oscillation.
\newblock \emph{Journal of Climate}, 14\penalty0 (1):\penalty0 98--115, 2001.

\bibitem[Timmermann et~al.(2003)Timmermann, Jin, and Abshagen]{timmermann2003nonlinear}
Axel Timmermann, Fei-Fei Jin, and Jan Abshagen.
\newblock A nonlinear theory for {E}l {N}i{\~n}o bursting.
\newblock \emph{Journal of the atmospheric sciences}, 60\penalty0 (1):\penalty0 152--165, 2003.

\bibitem[Roberts et~al.(2016)Roberts, Guckenheimer, Widiasih, Timmermann, and Jones]{roberts2016mixed}
Andrew Roberts, John Guckenheimer, Esther Widiasih, Axel Timmermann, and Christopher~KRT Jones.
\newblock Mixed-mode oscillations of {El}l {N}ino--{S}outhern oscillation.
\newblock \emph{Journal of the Atmospheric Sciences}, 73\penalty0 (4):\penalty0 1755--1766, 2016.

\bibitem[Timmermann and Jin(2002)]{timmermann2002nonlinear}
Axel Timmermann and Fei-Fei Jin.
\newblock A nonlinear mechanism for decadal {E}l {N}i{\~n}o amplitude changes.
\newblock \emph{Geophysical Research Letters}, 29\penalty0 (1):\penalty0 3--1, 2002.

\bibitem[Ren and Jin(2013)]{ren2013recharge}
Hong-Li Ren and Fei-Fei Jin.
\newblock Recharge oscillator mechanisms in two types of {ENSO}.
\newblock \emph{Journal of Climate}, 26\penalty0 (17):\penalty0 6506--6523, 2013.

\bibitem[Geng et~al.(2020)Geng, Cai, and Wu]{geng2020two}
Tao Geng, Wenju Cai, and Lixin Wu.
\newblock Two types of {ENSO} varying in tandem facilitated by nonlinear atmospheric convection.
\newblock \emph{Geophysical Research Letters}, 47\penalty0 (15):\penalty0 e2020GL088784, 2020.

\bibitem[Thual and Dewitte(2023)]{thual2023enso}
Sulian Thual and Boris Dewitte.
\newblock Enso complexity controlled by zonal shifts in the walker circulation.
\newblock \emph{Nature Geoscience}, 16\penalty0 (4):\penalty0 328--332, 2023.

\bibitem[Fang and Mu(2018)]{fang2018three}
Xiang-Hui Fang and Mu~Mu.
\newblock A three-region conceptual model for central {P}acific {E}l {N}i{\~n}o including zonal advective feedback.
\newblock \emph{Journal of Climate}, 31\penalty0 (13):\penalty0 4965--4979, 2018.

\bibitem[Chen et~al.(2022)Chen, Fang, and Yu]{chen2022multiscale}
Nan Chen, Xianghui Fang, and Jin-Yi Yu.
\newblock A multiscale model for {E}l {N}i{\~n}o complexity.
\newblock \emph{npj Climate and Atmospheric Science}, 5\penalty0 (1):\penalty0 1--13, 2022.

\bibitem[Jin et~al.(2007)Jin, Lin, Timmermann, and Zhao]{jin2007ensemble}
Fei-Fei Jin, L~Lin, A~Timmermann, and J~Zhao.
\newblock Ensemble-mean dynamics of the {ENSO} recharge oscillator under state-dependent stochastic forcing.
\newblock \emph{Geophysical Research Letters}, 34\penalty0 (3), 2007.

\bibitem[Chen and Majda(2017)]{chen2017simple}
Nan Chen and Andrew~J Majda.
\newblock Simple stochastic dynamical models capturing the statistical diversity of {E}l {N}i{\~n}o {S}outhern {O}scillation.
\newblock \emph{Proceedings of the National Academy of Sciences}, 114\penalty0 (7):\penalty0 1468--1473, 2017.

\bibitem[An and Jin(2004)]{an2004nonlinearity}
Soon-Il An and Fei-Fei Jin.
\newblock Nonlinearity and asymmetry of {ENSO}.
\newblock \emph{Journal of Climate}, 17\penalty0 (12):\penalty0 2399--2412, 2004.

\bibitem[Timmermann(2003)]{timmermann2003decadal}
Axel Timmermann.
\newblock Decadal enso amplitude modulations: A nonlinear paradigm.
\newblock \emph{Global and Planetary Change}, 37\penalty0 (1-2):\penalty0 135--156, 2003.

\bibitem[Jin et~al.(2020)Jin, Chen, Zhao, Hayashi, Karamperidou, Stuecker, Xie, and Geng]{jin2020simple}
Fei-Fei Jin, Han-Ching Chen, Sen Zhao, Michiya Hayashi, Christina Karamperidou, Malte~F Stuecker, Ruihuang Xie, and Licheng Geng.
\newblock Simple enso models.
\newblock \emph{El Ni{\~n}o Southern Oscillation in a changing climate}, pages 119--151, 2020.

\bibitem[Santosa and Symes(1986)]{santosa1986linear}
Fadil Santosa and William~W Symes.
\newblock Linear inversion of band-limited reflection seismograms.
\newblock \emph{SIAM journal on scientific and statistical computing}, 7\penalty0 (4):\penalty0 1307--1330, 1986.

\bibitem[Tibshirani(1996)]{tibshirani1996regression}
Robert Tibshirani.
\newblock Regression shrinkage and selection via the lasso.
\newblock \emph{Journal of the Royal Statistical Society Series B: Statistical Methodology}, 58\penalty0 (1):\penalty0 267--288, 1996.

\bibitem[AlMomani et~al.(2020)AlMomani, Sun, and Bollt]{almomani2020entropic}
Abd AlRahman~R AlMomani, Jie Sun, and Erik Bollt.
\newblock How entropic regression beats the outliers problem in nonlinear system identification.
\newblock \emph{Chaos: An Interdisciplinary Journal of Nonlinear Science}, 30\penalty0 (1), 2020.

\bibitem[Fish et~al.(2021)Fish, DeWitt, AlMomani, Laurienti, and Bollt]{fish2021entropic}
Jeremie Fish, Alexander DeWitt, Abd AlRahman~R AlMomani, Paul~J Laurienti, and Erik Bollt.
\newblock Entropic regression with neurologically motivated applications.
\newblock \emph{Chaos: An Interdisciplinary Journal of Nonlinear Science}, 31\penalty0 (11), 2021.

\bibitem[Kim et~al.(2017)Kim, Rogers, Sun, and Bollt]{kim2017causation}
Pileun Kim, Jonathan Rogers, Jie Sun, and Erik Bollt.
\newblock Causation entropy identifies sparsity structure for parameter estimation of dynamic systems.
\newblock \emph{Journal of Computational and Nonlinear Dynamics}, 12\penalty0 (1):\penalty0 011008, 2017.

\bibitem[AlMomani and Bollt(2020)]{almomani2020erfit}
Abd~AlRahman AlMomani and Erik Bollt.
\newblock Erfit: Entropic regression fit {MATLAB} package, for data-driven system identification of underlying dynamic equations.
\newblock \emph{arXiv preprint arXiv:2010.02411}, 2020.

\bibitem[Elinger(2021)]{elinger2021information}
Jared Elinger.
\newblock \emph{Information Theoretic Causality Measures For Parameter Estimation and System Identification.}
\newblock PhD thesis, Georgia Institute of Technology, Atlanta, GA, USA, 2021.

\bibitem[Elinger and Rogers(2021)]{elinger2021causation}
Jared Elinger and Jonathan Rogers.
\newblock Causation entropy method for covariate selection in dynamic models.
\newblock In \emph{2021 American Control Conference (ACC)}, pages 2842--2847. IEEE, 2021.

\bibitem[Chen(2020)]{chen2020learning}
Nan Chen.
\newblock Learning nonlinear turbulent dynamics from partial observations via analytically solvable conditional statistics.
\newblock \emph{Journal of Computational Physics}, page 109635, 2020.

\bibitem[Cover(1999)]{cover1999elements}
Thomas~M Cover.
\newblock \emph{Elements of information theory}.
\newblock John Wiley \& Sons, 1999.

\bibitem[Bellman(1961)]{bellman1961dynamic}
Richard~Ernest Bellman.
\newblock \emph{Dynamic programming treatment of the traveling salesman problem.}
\newblock RAND Corporation, 1961.

\bibitem[Majda and Chen(2018)]{majda2018model}
Andrew~J Majda and Nan Chen.
\newblock Model error, information barriers, state estimation and prediction in complex multiscale systems.
\newblock \emph{Entropy}, 20\penalty0 (9):\penalty0 644, 2018.

\bibitem[Tippett et~al.(2004)Tippett, Kleeman, and Tang]{tippett2004measuring}
Michael~K Tippett, Richard Kleeman, and Youmin Tang.
\newblock Measuring the potential utility of seasonal climate predictions.
\newblock \emph{Geophysical research letters}, 31\penalty0 (22), 2004.

\bibitem[Kleeman(2011)]{kleeman2011information}
Richard Kleeman.
\newblock Information theory and dynamical system predictability.
\newblock \emph{Entropy}, 13\penalty0 (3):\penalty0 612--649, 2011.

\bibitem[Branicki and Majda(2012)]{branicki2012quantifying}
Michal Branicki and Andrew~J Majda.
\newblock Quantifying uncertainty for predictions with model error in non-{G}aussian systems with intermittency.
\newblock \emph{Nonlinearity}, 25\penalty0 (9):\penalty0 2543, 2012.

\bibitem[Ying(2019)]{ying2019overview}
Xue Ying.
\newblock An overview of overfitting and its solutions.
\newblock In \emph{Journal of physics: Conference series}, volume 1168, page 022022. IOP Publishing, 2019.

\bibitem[Brunton et~al.(2016)Brunton, Proctor, and Kutz]{brunton2016discovering}
Steven~L Brunton, Joshua~L Proctor, and J~Nathan Kutz.
\newblock Discovering governing equations from data by sparse identification of nonlinear dynamical systems.
\newblock \emph{Proceedings of the national academy of sciences}, 113\penalty0 (15):\penalty0 3932--3937, 2016.

\bibitem[Chen and Zhang(2022)]{chen2022causality}
Nan Chen and Yinling Zhang.
\newblock A causality-based learning approach for discovering the underlying dynamics of complex systems from partial observations with stochastic parameterization.
\newblock \emph{arXiv preprint arXiv:2208.09104}, 2022.

\bibitem[Behringer and Xue(2004)]{behringer2004evaluation}
DW~Behringer and Yan Xue.
\newblock Evaluation of the global ocean data assimilation system at {NCEP}: {T}he {P}acific {O}cean.
\newblock In \emph{Proc. Eighth Symp. on Integrated Observing and Assimilation Systems for Atmosphere, Oceans, and Land Surface}, 2004.

\bibitem[Kalnay et~al.(1996)Kalnay, Kanamitsu, Kistler, Collins, Deaven, Gandin, Iredell, Saha, White, Woollen, et~al.]{kalnay1996ncep}
Eugenia Kalnay, Masao Kanamitsu, Robert Kistler, William Collins, Dennis Deaven, Lev Gandin, Mark Iredell, Suranjana Saha, Glenn White, John Woollen, et~al.
\newblock The {NCEP}/{NCAR} 40-year reanalysis project.
\newblock \emph{Bulletin of the American Meteorological Society}, 77\penalty0 (3):\penalty0 437--472, 1996.

\bibitem[Huang et~al.(2017)Huang, Thorne, Banzon, Boyer, Chepurin, Lawrimore, Menne, Smith, Vose, and Zhang]{huang2017extended}
Boyin Huang, Peter~W Thorne, Viva~F Banzon, Tim Boyer, Gennady Chepurin, Jay~H Lawrimore, Matthew~J Menne, Thomas~M Smith, Russell~S Vose, and Huai-Min Zhang.
\newblock Extended reconstructed sea surface temperature, version 5 (ersstv5): upgrades, validations, and intercomparisons.
\newblock \emph{Journal of Climate}, 30\penalty0 (20):\penalty0 8179--8205, 2017.

\bibitem[Thual et~al.(2016)Thual, Majda, Chen, and Stechmann]{thual2016simple}
Sulian Thual, Andrew~J Majda, Nan Chen, and Samuel~N Stechmann.
\newblock Simple stochastic model for {E}l {N}i{\~n}o with westerly wind bursts.
\newblock \emph{Proceedings of the National Academy of Sciences}, 113\penalty0 (37):\penalty0 10245--10250, 2016.

\bibitem[Palmer et~al.(2009)Palmer, Buizza, Doblas-Reyes, Jung, Leutbecher, Shutts, Steinheimer, and Weisheimer]{palmer2009stochastic}
Tim~N Palmer, Roberto Buizza, F~Doblas-Reyes, Thomas Jung, Martin Leutbecher, Glenn~J Shutts, Martin Steinheimer, and Antje Weisheimer.
\newblock Stochastic parametrization and model uncertainty.
\newblock \emph{ECMWF Technical Memorandum}, 598, 2009.

\bibitem[Chen et~al.(2015)Chen, Lian, Fu, Cane, Tang, Murtugudde, Song, Wu, and Zhou]{chen2015strong}
Dake Chen, Tao Lian, Congbin Fu, Mark~A Cane, Youmin Tang, Raghu Murtugudde, Xunshu Song, Qiaoyan Wu, and Lei Zhou.
\newblock Strong influence of westerly wind bursts on {E}l {N}i{\~n}o diversity.
\newblock \emph{Nature Geoscience}, 8\penalty0 (5):\penalty0 339--345, 2015.

\bibitem[Yang et~al.(2021)Yang, Majda, and Chen]{yang2021enso}
Qiu Yang, Andrew~J Majda, and Nan Chen.
\newblock {ENSO} diversity in a tropical stochastic skeleton model for the {MJO}, {E}l {N}i{\~n}o, and dynamic {W}alker circulation.
\newblock \emph{Journal of Climate}, pages 1--56, 2021.

\bibitem[Tziperman et~al.(1997)Tziperman, Zebiak, and Cane]{tziperman1997mechanisms}
Eli Tziperman, Stephen~E Zebiak, and Mark~A Cane.
\newblock Mechanisms of seasonal--{ENSO} interaction.
\newblock \emph{Journal of the Atmospheric Sciences}, 54\penalty0 (1):\penalty0 61--71, 1997.

\bibitem[Stein et~al.(2014)Stein, Timmermann, Schneider, Jin, and Stuecker]{stein2014enso}
Karl Stein, Axel Timmermann, Niklas Schneider, Fei-Fei Jin, and Malte~F Stuecker.
\newblock {ENSO} seasonal synchronization theory.
\newblock \emph{Journal of Climate}, 27\penalty0 (14):\penalty0 5285--5310, 2014.

\bibitem[Fang and Zheng(2021)]{fang2021effect}
Xiang-Hui Fang and Fei Zheng.
\newblock Effect of the air--sea coupled system change on the {ENSO} evolution from boreal spring.
\newblock \emph{Climate Dynamics}, pages 1--12, 2021.

\bibitem[Wang et~al.(2019)Wang, Luo, Yang, Sun, Cane, Cai, Yeh, and Liu]{wang2019historical}
Bin Wang, Xiao Luo, Young-Min Yang, Weiyi Sun, Mark~A Cane, Wenju Cai, Sang-Wook Yeh, and Jian Liu.
\newblock Historical change of {E}l {N}i{\~n}o properties sheds light on future changes of extreme {E}l {N}i{\~n}o.
\newblock \emph{Proceedings of the National Academy of Sciences}, 116\penalty0 (45):\penalty0 22512--22517, 2019.

\bibitem[Emerick and Reynolds(2013)]{emerick2013ensemble}
Alexandre~A Emerick and Albert~C Reynolds.
\newblock Ensemble smoother with multiple data assimilation.
\newblock \emph{Computers \& Geosciences}, 55:\penalty0 3--15, 2013.

\bibitem[Evensen(2003)]{evensen2003ensemble}
Geir Evensen.
\newblock The ensemble kalman filter: Theoretical formulation and practical implementation.
\newblock \emph{Ocean dynamics}, 53:\penalty0 343--367, 2003.

\bibitem[Thual et~al.(2019)Thual, Majda, and Chen]{thual2019statistical}
Sulian Thual, Andrew~J Majda, and Nan Chen.
\newblock Statistical occurrence and mechanisms of the 2014--2016 delayed super {E}l {N}i{\~n}o captured by a simple dynamical model.
\newblock \emph{Climate Dynamics}, 52\penalty0 (3-4):\penalty0 2351--2366, 2019.

\bibitem[Chen et~al.(2017)Chen, Li, Wang, and Wang]{chen2017formation}
Lin Chen, Tim Li, Bin Wang, and Lu~Wang.
\newblock Formation mechanism for 2015/16 super {E}l {N}i{\~n}o.
\newblock \emph{Scientific reports}, 7\penalty0 (1):\penalty0 1--10, 2017.

\bibitem[Hameed et~al.(2018)Hameed, Jin, and Thilakan]{hameed2018model}
Saji~N Hameed, Dachao Jin, and Vishnu Thilakan.
\newblock A model for super {E}l {N}i{\~n}os.
\newblock \emph{Nature Communications}, 9\penalty0 (1):\penalty0 2528, 2018.

\bibitem[Gubner(2006)]{gubner2006probability}
John~A Gubner.
\newblock \emph{Probability and random processes for electrical and computer engineers}.
\newblock Cambridge University Press, 2006.

\bibitem[Priestley(1981)]{priestley1981spectral}
Maurice~Bertram Priestley.
\newblock Spectral analysis and time series.
\newblock \emph{(No Title)}, 1981.

\bibitem[Kay(1988)]{kay1988modern}
Steven~M Kay.
\newblock \emph{Modern spectral estimation}.
\newblock Pearson Education India, 1988.

\bibitem[Percival and Walden(1993)]{percival1993spectral}
Donald~B Percival and Andrew~T Walden.
\newblock \emph{Spectral analysis for physical applications}.
\newblock cambridge university press, 1993.

\bibitem[Vialard et~al.(2001)Vialard, Menkes, Boulanger, Delecluse, Guilyardi, McPhaden, and Madec]{vialard2001model}
J{\'e}r{\^o}me Vialard, Christophe Menkes, Jean-Philippe Boulanger, Pascale Delecluse, Eric Guilyardi, Michael~J McPhaden, and Gurvan Madec.
\newblock A model study of oceanic mechanisms affecting equatorial pacific sea surface temperature during the 1997--98 el ni{\~n}o.
\newblock \emph{Journal of Physical Oceanography}, 31\penalty0 (7):\penalty0 1649--1675, 2001.

\bibitem[Jin and An(1999)]{jin1999thermocline}
Fei-Fei Jin and Soon-Il An.
\newblock Thermocline and zonal advective feedbacks within the equatorial ocean recharge oscillator model for enso.
\newblock \emph{Geophysical research letters}, 26\penalty0 (19):\penalty0 2989--2992, 1999.

\bibitem[Eisenman et~al.(2005)Eisenman, Yu, and Tziperman]{eisenman2005westerly}
Ian Eisenman, Lisan Yu, and Eli Tziperman.
\newblock Westerly wind bursts: {ENSO}'s tail rather than the dog?
\newblock \emph{Journal of Climate}, 18\penalty0 (24):\penalty0 5224--5238, 2005.

\bibitem[Levine and Jin(2017)]{levine2017simple}
Aaron~FZ Levine and Fei~Fei Jin.
\newblock A simple approach to quantifying the noise--enso interaction. part i: Deducing the state-dependency of the windstress forcing using monthly mean data.
\newblock \emph{Climate Dynamics}, 48\penalty0 (1-2):\penalty0 1--18, 2017.

\bibitem[Dieppois et~al.(2021)Dieppois, Capotondi, Pohl, Chun, Monerie, and Eden]{dieppois2021enso}
Bastien Dieppois, Antonietta Capotondi, Benjamin Pohl, Kwok~Pan Chun, Paul-Arthur Monerie, and Jonathan Eden.
\newblock {ENSO} diversity shows robust decadal variations that must be captured for accurate future projections.
\newblock \emph{Communications Earth \& Environment}, 2\penalty0 (1):\penalty0 1--13, 2021.

\bibitem[Yu and Kim(2013)]{yu2013identifying}
Jin-Yi Yu and Seon~Tae Kim.
\newblock Identifying the types of major {E}l {N}i{\~n}o events since 1870.
\newblock \emph{International Journal of Climatology}, 33\penalty0 (8):\penalty0 2105--2112, 2013.

\bibitem[Eyring et~al.(2016)Eyring, Bony, Meehl, Senior, Stevens, Stouffer, and Taylor]{eyring2016overview}
Veronika Eyring, Sandrine Bony, Gerald~A Meehl, Catherine~A Senior, Bjorn Stevens, Ronald~J Stouffer, and Karl~E Taylor.
\newblock Overview of the coupled model intercomparison project phase 6 (cmip6) experimental design and organization.
\newblock \emph{Geoscientific Model Development}, 9\penalty0 (5):\penalty0 1937--1958, 2016.

\bibitem[Brown et~al.(2008)Brown, Tudhope, Collins, and McGregor]{brown2008mid}
Josephine Brown, Alexander~W Tudhope, Matthew Collins, and Helen~V McGregor.
\newblock Mid-holocene enso: Issues in quantitative model-proxy data comparisons.
\newblock \emph{Paleoceanography}, 23\penalty0 (3), 2008.

\bibitem[Ford et~al.(2012)Ford, Ravelo, and Hovan]{ford2012deep}
Heather~L Ford, A~Christina Ravelo, and Steven Hovan.
\newblock A deep eastern equatorial pacific thermocline during the early pliocene warm period.
\newblock \emph{Earth and Planetary Science Letters}, 355:\penalty0 152--161, 2012.

\bibitem[Geng and Jin(2022)]{geng2022enso}
Licheng Geng and Fei-Fei Jin.
\newblock {ENSO} diversity simulated in a revised {C}ane-{Z}ebiak model.
\newblock \emph{Frontiers in Earth Science}, 10:\penalty0 899323, 2022.

\bibitem[Chen and Fang(2023)]{chen2023simple}
Nan Chen and Xianghui Fang.
\newblock A simple multiscale intermediate coupled stochastic model for el ni{\~n}o diversity and complexity.
\newblock \emph{Journal of Advances in Modeling Earth Systems}, 15\penalty0 (4):\penalty0 e2022MS003469, 2023.

\bibitem[Fang and Chen(2023)]{fang2023quantifying}
Xianghui Fang and Nan Chen.
\newblock Quantifying the predictability of enso complexity using a statistically accurate multiscale stochastic model and information theory.
\newblock \emph{Journal of Climate}, 36\penalty0 (8):\penalty0 2681--2702, 2023.

\bibitem[Chen(2023)]{chen2023stochastic}
Nan Chen.
\newblock \emph{Stochastic Methods for Modeling and Predicting Complex Dynamical Systems: Uncertainty Quantification, State Estimation, and Reduced-Order Models}.
\newblock Springer Nature, 2023.

\end{thebibliography}

\end{document}